\documentclass[12pt]{amsart}
\usepackage{div-alg-chars}
\usepackage{amsrefs,amscd}
\renewcommand\to{\oldto}
	\newcommand{\xrhardcode}[2]{\expandafter\def\csname#1\endcsname{#2}}
	\xrhardcode{exp-hyp:reduced}A
	\xrhardcode{exp-rem:when-hyps-hold}{1.2}
	\xrhardcode{exp-lem:complete-subfield}{2.2}
	\xrhardcode{exp-rem:actually-product}{2.1}
	\xrhardcode{exp-lem:torus-field-descent}{2.4}
	\xrhardcode{exp-lem:field-descent}{2.5}
	\xrhardcode{exp-lem:levi-descent}{2.6}
	\xrhardcode{exp-lem:domain-field-ascent}{2.7}
	\xrhardcode{exp-lem:unipotent}{2.8}
	\xrhardcode{exp-lem:stab-deep}{2.9}
	\xrhardcode{exp-rem:G-image}{3.1}
	\xrhardcode{exp-defn:compatibly-filtered}{3.3}
	\xrhardcode{exp-lem:ratl-maxl-torus}{3.2}
	\xrhardcode{exp-prop:compatibly-filtered-tame-rank}{3.4}
	\xrhardcode{exp-cor:compatibly-filtered-tame-rank}{3.6}
	\xrhardcode{exp-sec:normal}5
	\xrhardcode{exp-cor:good-comm}{6.2}
	\xrhardcode{exp-lem:x-depth}{6.4}
	\xrhardcode{exp-lem:split-in-center}{A.2}
\begin{document}
\begin{abstract}
The notion of a topological Jordan decomposition of a compact
element of a reductive $p$-adic group has proven useful in
many contexts.  In this paper, we generalise it to groups defined
over fairly general discretely valued fields and prove the
usual existence and uniqueness properties, as well as an
analogue of a fixed-point result of Prasad and Yu.
\end{abstract}

\title{Topological Jordan decompositions}
\subjclass[2000]{Primary 20G25, 22E35.  Secondary 22D05.}
\author{Loren Spice}
\address{The University of Michigan \\
Ann Arbor, MI 48109-1043}
\email{lspice@umich.edu}
\thanks{The author was supported by
National Science Foundation Postdoctoral Fellowship award
DMS-0503107.}
\date\today
\maketitle

\addtocounter{section}{-1}

\section{Introduction}

In \cite{kazhdan:lifting}, Kazhdan defines the
notions of \ff-semisimplicity and \ff-unipotence of an
element of $\GL_n(F_0)$, where $F$ is a discretely valued
locally compact field with ring of integers $F_0$ and residue field
\ff (see the definition on p.~226 of \emph{loc.\ cit.}).
An arbitrary element of $\GL_n(F_0)$ can be decomposed
as a commuting product of an \ff-semisimple and an
\ff-unipotent element (see Lemma 2 on p.~226 of \emph{loc.\
cit.}).
Furthermore, stably conjugate \ff-semisimple elements are
actually $\GL_n(F_0)$-conjugate
(see Lemma 3 on p.~226 of \emph{loc.\ cit.}, where the
result is proven for rationally conjugate elements, and
Lemma 13.1 of \cite{hales:simple-defn}).

Kazhdan uses this last result in his calculation of the
$\varepsilon$-twisted orbital integral $I_\ell(f)$ (see
Theorem 1 on p.~224 of \cite{kazhdan:lifting}, and the definition
immediately preceding it).
A detailed exposition appears in
\cite{henniart:ll-cyclic}*{\S 5}; see especially
\S 5.6 of \emph{loc.\ cit.}
An analogous result is used by Waldspurger in his
computation of Shalika germs for $\GL(n)$ (see
\cite{waldspurger:shalika}*{\S 5}).

In \cite{hales:simple-defn}, Hales defines absolute
semisimplicity and topological unipotence (the analogues of
\ff-semisimplicity and \ff-unipotence) for elements of
unramified groups, and shows that every strongly compact
element can be decomposed as a commuting product of an
absolutely semisimple and a topologically unipotent element
(the topological Jordan decomposition).
He then defines transfer factors on unramified groups, and
shows that the transfer factor at a strongly compact element may
be expressed in terms of the transfer factor for the centraliser
of its absolutely semisimple part, evaluated at its
topologically unipotent part (see Theorem 10.18 and Lemma
13.2 of \emph{loc.\ cit.}).
This suggests that the topological Jordan decomposition
is important for the fundamental lemma.
Indeed, in \cite{flicker:sym-square}, Flicker uses a twisted
analogue of the decomposition to prove a special case of the
fundamental lemma (see the theorem on p.~509 of \emph{loc.\ cit.}).

The topological Jordan decomposition is also useful in
character computations.
Recall that the character of a Deligne--Lusztig
representation of a finite group of Lie type is expressed by a reduction
formula in terms of the (ordinary) Jordan decomposition (see
Theorem 4.2 of \cite{deligne-lusztig:finite}).
The topological Jordan decomposition plays the same r\^ole
for the characters of depth-zero supercuspidal
representations of $p$-adic groups arising via compact
induction from representations which are inflations of
Deligne--Lusztig representations of reductive quotients
(see Lemma 10.0.4 of \cite{debacker-reeder:depth-zero-sc}).

In
\cite{adler-spice:good-expansions}*{\S\xref{exp-sec:normal}},
as preparation for the (positive-depth) character computations of
\cite{adler-spice:explicit-chars},
Adler and the author define the notion of a \emph{normal
approximation} of an element of a reductive $p$-adic group,
a refinement of the topological Jordan
decomposition.  However, for the results of that paper, one
needs notions of absolute semisimplicity
and topological unipotence
that make sense over a discretely valued field $F$ which
is not necessarily locally compact.

In this paper, we offer two generalisations of these
notions, the first an abstract one adapted to
profinite groups, and the second adapted to the setting in
which we are most interested, of reductive groups over
discretely valued fields $F$ as above.  We prove the usual
existence (Propositions \ref{prop:intrinsic-tu} and
\ref{prop:tJd=cpct}) and uniqueness
(Propositions \ref{prop:unique-top-p-Jordan} and
\ref{prop:unique-top-F-Jordan}) results for these
decompositions.
In the familiar case where $F$ is locally compact, these
definitions have significant overlap; see, for example,
Lemmata \ref{lem:top-unip}, \ref{lem:abs-ss}, and
\ref{lem:top-Jordan}.
Our main result, Theorem \ref{thm:top-F-Jordan-projects}, is
a strong existence result which is the analogue of item (7) in the
list of properties of topological Jordan decompositions
given in \cite{hales:simple-defn}*{\S 3}.
An important ingredient in its proof is an analogue of a
fixed-point result of Prasad and Yu (see Proposition
\ref{prop:x-depth}).

This work has its origins in discussions with Jeffrey Adler
during the writing of \cite{adler-spice:good-expansions}.  I
thank him for careful reading and many useful suggestions.
I also thank Stephen DeBacker and Gopal Prasad for helpful
conversations.
Finally, I thank the referee for many valuable comments which
allowed me to improve several proofs and correct some confusing
typos.

\section{Abstract groups}
\label{sec:profinite}

Fix a prime $p$, a Hausdorff topological group $G$, and a
closed normal subgroup $N$.
Note that $G/N$ is also Hausdorff.
For $g, \gamma \in G$, we define
$\lsup g\gamma := g\gamma g\inv$.

\begin{dn}
An element or subgroup of $G$ is \emph{compact modulo $N$}
if its image in $G/N$ belongs to a compact subgroup.
If $N$ is the trivial subgroup, we shall omit ``modulo $N$''.
\end{dn}

\begin{dn}
\label{defn:ILPF}
The group $G$ is \emph{ind-locally-compact}
(respectively, \emph{ind-locally-profinite};
respectively, \emph{ind-locally-pro-$p$})
if it is an inductive limit of a directed system of
locally compact
(respectively, locally profinite; respectively, locally pro-$p$)
groups.
\end{dn}

The main example we will have in mind of an
ind-locally-compact group is
the set of $F$-rational points of a linear algebraic $F$-group \bG,
where $F$ is an algebraic extension of a locally compact
field (see Remark \ref{rem:ILPF}).
Another example is
$\GL_\infty(F) := \varinjlim \GL_n(F)$,
where $F$ is a locally compact field;
the limit is taken over positive integers $n$;
and, for a given $n$,
the map $\GL_n(F) \to \GL_{n + 1}(F)$
comes from the natural embedding of $\GL_n$
into the Levi subgroup $\GL_n \times \GL_1$
of $\GL_{n + 1}$.

\begin{dn}
\label{defn:top-p-ss-unip}
An element $\gamma \in G$
is \emph{absolutely $p$-semisimple}
if it has finite, coprime-to-$p$ order.
It is \emph{topologically $p$-unipotent} if
$\lim_{n \to \infty} \gamma^{p^n} = 1$.
If the projection of $\gamma$ to $G/N$ is absolutely
$p$-semisimple (respectively, topologically $p$-unipotent),
then we will say that $\gamma$ is \emph{absolutely $p$-semisimple
modulo $N$} (respectively, \emph{topologically $p$-unipotent modulo
$N$}).
\end{dn}

\begin{rk}
\label{rem:power-of-p-ss-unip}
Any power of an absolutely $p$-semisimple modulo $N$ (respectively,
topologically $p$-unipotent modulo $N$) element is absolutely
$p$-semisimple modulo $N$ (respectively, topologically $p$-unipotent
modulo $N$).
A $p$-power root of a topologically $p$-unipotent modulo $N$
element is again topologically $p$-unipotent modulo $N$.
\end{rk}

\begin{rk}
\label{rem:top-p-ss-unip-cpct}
Clearly, an absolutely $p$-semisimple element is compact.
Suppose that $\gamma \in G$ is topologically $p$-unipotent
and $G$ is ind-locally-compact.
Then $\gamma$ belongs to a locally compact subgroup of $G$,
so, since $\gamma^{p^n} \to 1$, there is some $n \in \Z_{> 0}$
such that $\gamma^{p^n}$ belongs to a compact subgroup of
$G$.  Thus $\gamma$ is compact.
\end{rk}

\begin{dn}
\label{defn:top-p-Jordan}
A \emph{topological $p$-Jordan decomposition modulo $N$}
of an element $\gamma \in G$
is a pair of commuting elements $(\gamma\tsemi, \gamma\tunip)$
of $G$ such that
\begin{itemize}
\item $\gamma = \gamma\tsemi\gamma\tunip$,
\item $\gamma\tsemi$ is absolutely $p$-semisimple modulo $N$,
and
\item $\gamma\tunip$ is topologically $p$-unipotent modulo $N$.
\end{itemize}
We will sometimes just say that
$\gamma = \gamma\tsemi\gamma\tunip$ is a topological $p$-Jordan
decomposition modulo $N$.
If $N$ is the trivial subgroup, we will omit ``modulo $N$''.
\end{dn}

In the statement of the following result, recall that $p$ is
\emph{fixed}.  It is certainly possible for an element to have
distinct $p$- and $\ell$-decompositions for $\ell$ a
prime distinct from $p$
(although it is an easy consequence of Remark
\ref{rem:intrinsic-tu} that, if $G$ is
ind-locally-pro-$p$, then this happens only for finite-order
elements).

\begin{pn}
\label{prop:unique-top-p-Jordan}
Suppose that $\gamma \in G$ has a topological $p$-Jordan
decomposition $\gamma = \gamma\tsemi\gamma\tunip$.
\begin{inc_enumerate}
\item
If $\gamma = \gamma\tsemi'\gamma\tunip'$ is a topological
$p$-Jordan decomposition, then
$\gamma\tsemi = \gamma\tsemi'$ and
$\gamma\tunip = \gamma\tunip'$.
\item\label{prop:unique-top-p-Jordan_closure}
The closure of the group generated by $\gamma$ contains
$\gamma\tsemi$ and $\gamma\tunip$.
\item
If $G'$ is another Hausdorff topological group,
and $f : G \to G'$ is a continuous homomorphism, then
$f(\gamma) = f(\gamma\tsemi)f(\gamma\tunip)$ is a
topological $p$-Jordan decomposition.
\item
For $g \in G$, we have that
$\lsup g\gamma = (\lsup g\gamma\tsemi)(\lsup g\gamma\tunip)$
is a topological $p$-Jordan decomposition.
\end{inc_enumerate}
\end{pn}

A special case of the above was introduced
in \cite{hales:simple-defn}*{\S 3}
(especially items (3) and (4) of the list there), and
Lemma 2 on \cite{kazhdan:lifting}*{p.~226}.
We omit the (straightforward) proof.


\begin{pn}
\label{prop:intrinsic-tu}
If $G$ is ind-locally-profinite, then
an element $\gamma\in G$ is topologically $p$-unipotent
if and only if
it belongs to a pro-$p$ subgroup of $G$.
If $G$ is ind-locally-pro-$p$, then an
element $\gamma \in G$ has a topological $p$-Jordan
decomposition if and only if it is compact.
\end{pn}

\begin{proof}
Suppose that $G$ is ind-locally-profinite.
Note that $G$ is ind-locally-compact and totally
disconnected.
Denote by $K$ the closure of the subgroup of $G$ generated
by $\gamma$.  Then $K$ is also totally disconnected.
By Proposition I.1.1.0 of \cite{serre:galois}, if $K$ is
compact, then it is profinite.
In this case,
since $K$ is Abelian,
by Proposition I.1.4.3 of \cite{serre:galois}
it has a unique Sylow pro-$\ell$ subgroup for each
prime $\ell$.
Let $K_p$ be its Sylow pro-$p$ subgroup,
and $K_{p'}$ the direct product of its Sylow pro-$\ell$
subgroups, taken over all $\ell \ne p$.
Then $K_p$ and $K_{p'}$ are profinite Abelian groups, of
$p$-power and prime-to-$p$ order, respectively.
It is an easy consequence of Proposition I.1.4.4(b) of
\cite{serre:galois} that $K = K_p \times K_{p'}$.
Write $\gamma = \gamma_p\gamma_{p'}$, with
$\gamma_p \in K_p$ and $\gamma_{p'} \in K_{p'}$.

The `if' direction of the first statement is obvious.

For the `only if' direction of the first statement,
suppose that $\gamma$ is topologically $p$-unipotent.
By Remark \ref{rem:top-p-ss-unip-cpct}, $K$ is compact.
By the first paragraph of the proof, it is profinite.
For any open subgroup
$U'$ of $K_{p'}$, we have that the map $g \mapsto g^p$ is an
isomorphism on $K_{p'}/U'$.
(Here, we have used commutativity of $K_{p'}$.)
Since $\gamma_{p'}^{p^n} \in U'$ for some $n \in \Z_{> 0}$, we
have that $\gamma_{p'} \in U'$.
Since $U'$ was arbitrary and $K_{p'}$ is Hausdorff, we have that
$\gamma_{p'} = 1$; i.e., $\gamma = \gamma_p \in K_p$.
Since $K_p$ is closed, in fact $K = K_p$; i.e., $K$ is a
pro-$p$ group.

Now suppose that $G$ is ind-locally-pro-$p$.
The `only if' direction of the second statement follows from
Remark \ref{rem:top-p-ss-unip-cpct}.

For the `if' direction of the second statement, suppose that
$\gamma$ (equivalently, $K$) is compact.
By the first paragraph of the proof, $K$ is profinite.
Then the intersection of $K_{p'}$ with an open ind-pro-$p$
subgroup of $H$ is an open ind-pro-$p$ subgroup of $K_{p'}$.
However, $K_{p'}$ contains no non-trivial pro-$p$ subgroup, so
this intersection is the trivial subgroup of $K_{p'}$.
Thus $K_{p'}$ is discrete, so finite.
Put $\gamma\tsemi = \gamma_{p'}$
and $\gamma\tunip = \gamma_p$.
By the first statement of the lemma, $\gamma\tunip$ is topologically
$p$-unipotent.
Since $K_{p'}$ is finite, $\gamma\tsemi$ is absolutely
$p$-semisimple.
\end{proof}

\begin{rk}
\label{rem:intrinsic-tu}
We isolate from the preceding proof a more refined, but
technical, version of Proposition
\ref{prop:intrinsic-tu}.  Suppose that $G$ is
ind-locally-profinite, $\gamma$ is a compact element of $G$, and
$K$ is the closure in $G$ of the group generated by
$\gamma$.  Then $K$ is profinite, and we may write
$K = K_p \times K_{p'}$, where $K_p$ and $K_{p'}$ are
profinite groups of $p$-power and prime-to-$p$ order,
respectively.  Write
$\gamma = \gamma_p\gamma_{p'}$, with
$\gamma_p \in K_p$ and $\gamma_{p'} \in K_{p'}$.
Then
\begin{itemize}
\item $\gamma$ is topologically $p$-unipotent if and only if
$\gamma_{p'} = 1$, in which case $K_{p'} = \sset 1$.
\item $\gamma$ has a topological $p$-Jordan decomposition if
and only if $\gamma_{p'}$ has finite order, in which case
$K_{p'}$ is finite and
we may take the topologically $p$-semisimple and
topologically $p$-unipotent parts of $\gamma$ to be
$\gamma_{p'}$ and $\gamma_p$, respectively.
\item If $G$ is ind-locally-pro-$p$, then $K_{p'}$
is finite, so $\gamma$ has a topological $p$-Jordan
decomposition.
\end{itemize}
\end{rk}

\section{Algebraic groups}

\subsection{Unipotent elements}

The adjective ``unipotent'' has sometimes carried several
meanings (see \cite{adler-debacker:bt-lie}*{\S 3.7.1}).
We begin by defining our notion of unipotence, then give a
general result, essentially due to Kempf,
relating the different meanings.

\begin{dn}
If $F$ is a field, \bG is a linear algebraic $F$-group,
and $\gamma \in \bG(F)$, then $\gamma$ is \emph{unipotent}
if there is
an embedding $\bG \hookrightarrow \GL_n$, for some
$n \in \Z_{\ge 0}$, such that the image of $\gamma$ is
an upper triangular matrix, with $1$s on the diagonal.
\end{dn}

\begin{lm}
\label{lem:equiv-unip}
If $F$ is a field, \bG is a connected reductive $F$-group,
and $\gamma \in \bG(F)$ is unipotent, then there are a
finite separable extension $E/F$ and a unipotent radical \bU
of a parabolic $E$-subgroup of \bG such that
$\gamma \in \bU(E)$.
If $E$ is equipped with a topology making it a non-discrete
Hausdorff topological field,
then there is a one-parameter subgroup $\lambda$ of \bG,
defined over $E$,
such that $\lim_{t \to 0} \lsup{\lambda(t)}\gamma = 1$ in the
$E$-analytic topology on $\bG(E)$.
\end{lm}

Since \bG embeds as a closed subset of the affine space
$\mb A^N$ for some $N \in \Z_{\ge 0}$, we may regard $\bG(E)$
as a subset of $E^N$.  By definition, the $E$-analytic topology
on $\bG(E)$ is just the subspace topology.
This topology is finer than the Zariski topology on $\bG(E)$.
By \cite{conrad:finite-class-number}*{Appendix B}
(especially Theorem B.1),
it is independent of the choice of embedding.

\begin{proof}
By Lemma 3.1 and Theorem 3.4 of \cite{kempf:instability}, there is a
one-parameter subgroup $\lambda$ of \bG such that
$\lim_{t \to 0} \lsup{\lambda(t)}\gamma = 1$
in the Zariski topology.
(Indeed, in the notation of Theorem 3.4(c) of \emph{loc.\ cit.}, any
$\lambda \in \Delta_{S = \sset 1, x = \gamma}$ will do.)
There is a finite separable extension $E/F$ such that
\bG is $E$-split and $\lambda$ is defined over $E$.
Then $\gamma \in \bU(\lambda)(E)$,
in the notation of \cite{kempf:instability}*{p.~305}.
By Theorem 13.4.2(i) and Lemma 15.1.2(ii) of
\cite{springer:lag}, $\bU(\lambda)$ is the unipotent radical
of a parabolic $E$-subgroup of \bG.

Now let \bT be an $E$-split maximal torus containing the
image of $\lambda$, and $\Phi(\bG, \bT)$ the root system of \bT
in \bG.  Put
$\Phi_\lambda^+
= \set{\alpha \in \Phi(\bG, \bT)}
	{\langle\alpha, \lambda\rangle > 0}$.
By Proposition 14.4(2)(a) of \cite{borel:linear},
$\bU(\lambda)$ is (as a variety) the Cartesian product of
the root subgroups $\bU_\alpha$ of \bG associated to roots
$\alpha \in \Phi_\lambda^+$.
By Theorem 18.7 of \cite{borel:linear}, there are $E$-isomorphisms
$\mexp_\alpha : \Add \to \bU_\alpha$
for $\alpha \in \Phi_\lambda^+$ such that
$\lsup \tau\mexp_\alpha(s) = \mexp_\alpha(\alpha(\tau)s)$
for $s \in E$ and $\tau \in \bT(E)$.
There are elements
$s_\alpha \in E$ (for $\alpha \in \Phi_\lambda^+$) such that
$\gamma = \prod_{\alpha \in \Phi_\lambda^+}
	\mexp_\alpha(s_\alpha)$.
Now suppose that $E$ is equipped with a topology making it a
non-discrete Hausdorff
topological field.
By Theorem B.1 of \cite{conrad:finite-class-number},
the maps
$\mexp_\alpha : E \to \bU_\alpha(E) \subseteq \bG(E)$ are
continuous for the $E$-analytic topology on $\bG(E)$,
so
$$
\lim_{t \to 0} \lsup{\lambda(t)}\gamma
= \prod_{\alpha \in \Phi_\lambda^+} \lim_{t \to 0}
	\lsup{\lambda(t)}\mexp_\alpha(s_\alpha)
= \prod_{\alpha \in \Phi_\lambda^+}
	\mexp_\alpha\Bigl(\lim_{t \to 0}
		t^{\langle\alpha, \lambda\rangle}s_\alpha
	\Bigr)
= \prod_{\alpha \in \Phi_\lambda^+} \mexp_\alpha(0)
= 1
$$
in the $E$-analytic topology on $\bG(E)$, as desired.
\end{proof}

\subsection{Algebraic groups:  basic definitions and notation}

Let
\begin{itemize}
\item $F$ be a field, with non-trivial discrete valuation
$\ord$, that is an algebraic extension of a complete field
with perfect residue field,
\item \ol F an algebraic closure of $F$,
\item $F\unram/F$ the maximal unramified subextension of
$\ol F/F$,
\item $F\tame/F$ the maximal tame subextension of $\ol F/F$,
\item $F_0$ the ring of integers of $F$,
\item $F_{0+}$ the maximal ideal of $F_0$\,,
\item $F\cross_0 = F_0 \smallsetminus F_{0+}$\,,
\item $F\cross_{0+} = 1 + F_{0+}$\,,
\item \ff the residue field $F_0/F_{0+}$\,,
\item $p = \chr \ff$,
\item
$(F\cross)^{p^\infty}
= \bigcap_{n = 0}^\infty
	(F\cross)^{p^n}$,
\item \bG a connected reductive $F$-group,
\item \mo N a closed normal $F$-subgroup of \bG,
and
\item \tbG the quotient $\bG/\mo N$.
\end{itemize}
(Many of our results apply to any Henselian field
with perfect residue field; but we restrict our attention slightly
so that we do
not have to re-prove the results of
\cite{adler-spice:good-expansions}*{\S 2} in this
generality.)
We denote by $Z(\bG)$ the centre of \bG;
by $\bX^*(\bG)$ and $\bX_*(\bG)$ the
characters and cocharacters, respectively, of \bG; and by
$\bX^*_F(\bG)$ and $\bX_*^F(\bG)$ those characters and
cocharacters, respectively, defined over $F$.
If necessary, we will write $\ff_F$ in place of \ff to
indicate the dependence on the field $F$.

We will assume without further mention that any algebraic
extension of $F$ is contained in \ol F.
If $E/F$ is such an extension, then
we denote again by $\ord$ the unique extension of $\ord$ to
a (not necessarily discrete) valuation on $E$;
and by $E_0$\,, etc., the analogues for $E$ of $F_0$\,, etc.,
above.

We will write $G = \bG(F)$ and $N = \mo N(F)$, and similarly
for other $F$-groups.

\begin{dn}
If $\gamma \in G$ is semisimple, then the
\emph{character values} of $\gamma$ (in \bG) are the
elements of the set
$\set{\chi(\gamma)}{\chi \in \bX^*(\bT)}$,
where \bT is any maximal torus in \bG containing $\gamma$.
\end{dn}

\begin{dn}
\label{defn:tame}
An element $\gamma \in G$ is \emph{$F$-tame} if there exists
an $F\tame$-split torus
(equivalently, by Lemmata \xref{exp-lem:ratl-maxl-torus} and
\xref{exp-lem:split-in-center} of
\cite{adler-spice:good-expansions}, an $F\tame$-split
$F$-torus)
\bS in \bG such that $\gamma \in \bS(F\tame)$.
\end{dn}

\begin{dn}
Let $\BB(\bG, F)$ be the (enlarged) Bruhat--Tits building of
\bG over $F$ and, for $x \in \BB(\bG, F)$, let
$G_x$ and $G_x^+$ be the parahoric subgroup associated to
$x$ and its pro-unipotent radical, respectively.
(In general, the parahoric subgroup may be strictly smaller
than the stabiliser of $x$ (but see Lemma
\ref{lem:tame-para-ascent}).
In the language of Proposition
4.6.28(i) of \cite{bruhat-tits:reductive-groups-2}, it is the
\emph{fixateur connexe} of the facet containing $x$.)
Let $\ms G_x$ be the (not necessarily connected) $\ff_F$-group such that
$\ms G_x(\ff_{\wtilde F})
= \stab_{\bG(\wtilde F)}(x)/\bG(\wtilde F)_x^+$
for all unramified extensions $\wtilde F/F$.
If necessary, we will write $\ms G_x^F$ in place of $\ms G_x$
to indicate the dependence on the field $F$.
Put
$G_0 = \bigcup_{x \in \BB(\bG, F)} G_x$ and
$G_{0+} = \bigcup_{x \in \BB(\bG, F)} G_x^+$\,.
\end{dn}

\begin{rk}
We have
$\ms G_x\conn(\ff_{\wtilde F}) = \bG(\wtilde F)_x/\bG(\wtilde F)_x^+$
for all unramified extensions $\wtilde F/F$.
\end{rk}

\begin{dn}
An element or subgroup of $G$ is
\emph{bounded} if its orbits in $\BB(\bG, F)$ are bounded
(in the sense of metric spaces).
An element or subgroup of $G$ is \emph{bounded modulo \mo N{}}
if its image in \tG is bounded.
If $\bG = \bT$ is a torus, then denote by $T\subb$ the
maximal bounded subgroup of $T$.
\end{dn}

\begin{rk}
If \tbG is semisimple, then the building
$\BB(\tbG, F)$ is canonical.  In
general, we ``canonify'' it as in
\cite{tits:corvallis}*{\S\S 1.2 and 2.1}.
Since we will be concerned almost exclusively with the case
$\mo N = Z(\bG)\conn$, this ``canonification'' will not
usually be necessary.
\end{rk}

\begin{rk}
\label{rem:equiv-bdd}
Consider a bounded element or subgroup of $G$ and a
non-empty, closed, convex, $G$-stable subset \mc S of $\BB(\bG, F)$.
By Proposition 3.2.4 of
\cite{bruhat-tits:reductive-groups-1}, the element or
subgroup fixes a point \ox of the image of \mc S in
the reduced building $\rBB(\bG, F)$; hence, by boundedness, actually
fixes any lift to \mc S of \ox.
On the other hand, since $G$ acts on $\BB(\bG, F)$ by isometries,
an element or subgroup of $G$
which fixes a point $x \in \BB(\bG, F)$ is bounded.
\end{rk}

\begin{rk}
\label{rem:bdd=cpct}
If $F$ is locally compact,
then a subgroup of $G$ is bounded modulo \mo N if and only if its
closure is compact modulo $N$.
If $F$ is an algebraic extension of a locally compact field,
then an element of $G$ is bounded modulo \mo N if and only
if it belongs to a compact modulo $N$ subgroup of $G$.
Indeed, the `if' direction is obvious.
For the `only if' direction, suppose that $\gamma \in G$ is
bounded modulo \mo N.  Then, by Lemma
\xref{exp-lem:complete-subfield} of
\cite{adler-spice:good-expansions}, there is some locally
compact subfield $F'$ of $F$ such that \bG, \mo N, and
$\gamma$ are all defined over $F'$.
Thus $\gamma$ is contained in a compact modulo $\mo N(F')$
subgroup of $\bG(F')$, hence \emph{a fortiori} a compact
modulo $N$ subgroup of $G$.
\end{rk}

\begin{rk}
\label{rem:torus-filt}
If $\bG = \bT$ is a torus, then
$T_0 = T_x$ and $T_{0+} = T_x^+$ for any
$x \in \BB(\bT, F)$.
Concretely,
$T\subb$ is the group of elements of $T$ whose character
values lie in $E\cross_0$ (by Lemme 4.2.19 of
\cite{bruhat-tits:reductive-groups-2})
and
$T_x^+$ is the group
of elements of $T_x$ whose character values lie in
$E\cross_{0+}$\,, where $E/F$ is the splitting
field of \bT.
If \bT is $F$-split, then
$T_x = T\subb$\,; so $T_x^+$ is the group of
elements of $T$ whose character values lie in
$F\cross_{0+}$\,.
\end{rk}

\begin{rk}
\label{rem:ILPF}
If $F$ is an algebraic extension of a locally compact field,
then $p > 0$ and, by Lemma \xref{exp-lem:complete-subfield} of
\cite{adler-spice:good-expansions}, we have
$G = \varinjlim \bG(F')$, the limit taken over all locally
compact subfields $F'$ of $F$ over which \bG is defined.
For such a subfield, $\bG(F')_x^+$ is an open pro-$p$
subgroup of $\bG(F')$ (for any $x \in \BB(\bG, F')$).  Thus,
$G$ is ind-locally-pro-$p$.
\end{rk}

\begin{rk}
\label{rem:equiv-unip}
Suppose that $\gamma \in G$ is unipotent.  By Lemma
\ref{lem:equiv-unip}, there are a finite separable extension
$E/F$ (which we may take, by passing to a further finite
separable extension
if necessary, to be a splitting field for \bG)
and a one-parameter subgroup
$\lambda \in \bX_*^E(\bG)$ such that 
$\lim_{t \to 0} \lsup{\lambda(t)}\gamma = 1$.
Since $\bG(E)_{0+}$ is a neighbourhood of
$1$, there is an element $t \in E\cross$ such that
$\lsup{\lambda(t)}\gamma \in \bG(E)_{0+}$\,.
Then
$\gamma \in \lsup{\lambda(t)\inv}\bG(E)_{0+}
	= \bG(E)_{0+}$\,,
so $\gamma$ is topologically $F$-unipotent.
\end{rk}

\subsection{A lifting of $\smash{\ol\ff}\cross$}
\label{subsec:coefficient}

In this section, we will define a
$\Gal(F\unram/F)$-stable subgroup
$\mc F(F)$ of $(F\unram)\cross$ such that
the map
$\mc F(F) \to \smash{\ol\ff}\cross$ is an isomorphism.
Let $F'$ be a complete subfield of $F$ such that $F/F'$ is
unramified.  (Such a subfield exists, by Lemma
\xref{exp-lem:complete-subfield} of
\cite{adler-spice:good-expansions}.)

If $p > 0$, then put
$\mc F(F)
:= \bigcup_{\text{$L/F'$ finite unramified}}
	(L\cross)^{p^\infty}$.
We have that
$\mc F(F) \cap (F\unram)\cross_{0+} = \sset 1$,
so that the map $\mc F(F) \to \smash{\ol\ff}\cross$
is injective.
By Lemma 7 of \cite{cohen:complete-local-rings},
since \ff is perfect, $(L\cross)^{p^\infty}$
maps onto $\ff_L\cross$ for every finite unramified
extension $L/F'$; so the map
$\mc F(F) \to \smash{\ol\ff}\cross$ is also surjective,
hence again an isomorphism.
Note that the $p^n$th power map on $\mc F(F)$ is also an
isomorphism for all $n \in \Z_{\ge 0}$.
Note that, if $E/F'$ is an arbitrary finite extension with
maximal unramified subextension $L/F'$, then
$(E\cross)^{p^\infty}$ contains $(L\cross)^{p^\infty}$,
and both map isomorphically onto $\ff_E = \ff_L$;
so, in fact, $(E\cross)^{p^\infty} = (L\cross)^{p^\infty}$,
and we could take the union defining $\mc F(F)$ over
\emph{all} finite extensions $E/F'$.

The definition of $\mc F(F)$ is slightly more complicated if
$p = \chr F$.
Let $\ff_F'$ be a subfield of $F'_0$ satisfying the
following property.
\begin{itemize}
\item[($\text{\textbf{CF}}_F$)]  The restriction to $\ff_F'$ of the
natural map $F'_0 \to \ff_F$ is an isomorphism onto $\ff_F$.
\end{itemize}
By Theorem 9 of \cite{cohen:complete-local-rings},
$\ff_F'$ exists.

If $L/F$ is a finite unramified extension,
say of degree $n$,
then $\ff_L/\ff_F$ is separable, so there
exists a primitive element $\ol\theta$ for
$\ff_L/\ff_F$, say with minimal polynomial $\ol m(x)$ over $\ff_F$.
Since $\ff_L/\ff_F$ is separable, so is $\ol m(x)$.
Let $m(x)$ be the unique preimage in $\ff_F'[x]$ of $\ol m(x)$,
and $\theta$ the unique root of $m(x)$ lifting $\ol\theta$.
Note that $F'[\theta]$ is a finite unramified extension of $F'$,
hence complete.

Suppose that there exists a subfield $\ff_L'$ of $L_0$,
containing $\ff_F'$, with property ($\text{\textbf{CF}}_L$).
Then $\indx{\ff_L'}{\ff_F'} = n$, and
there is a lift $\theta'$ in $\ff_L'$ of
$\ol\theta$, say with minimal polynomial $m'(x)$ over $\ff_F'$.
Then $\deg m'(x) \le \indx{\ff_L'}{\ff_F'} = n$ and
$\ol\theta$ is a root of the image in $\ff_F[x]$ of
$m'(x)$, so $m'(x)$ is the preimage in
$\ff_F'[x]$ of $\ol m(x)$; that is, $m'(x) = m(x)$.
Thus $\theta' = \theta$, so $\ff_L' = \ff_F'[\theta]$.
Note that $\ff_L'$ lies in the complete field $F'[\theta]$.

Since $\ff_F'[\theta]$ clearly has property
($\text{\textbf{CF}}_L$), we have shown that it is the
unique subfield of $L_0$ containing $\ff_F'$ with this property.
In particular, if $L/F$ is Galois (which is not automatic,
since we have not assumed that $\ff_F$ is finite), then
$\ff_L' \subseteq L$ is $\Gal(L/F)$-,
hence $\Gal(F\unram/F)$-, stable.
Put
$\mc F(F) := \bigcup_{\text{$L/F$ finite unramified}}
	(\ff_L')\cross$.
By Theorem 10(b) of \cite{cohen:complete-local-rings}, if
$p > 0$ (in addition to $p = \chr F$), then this definition
coincides with the one given above.

It is clear that $\mc F(F) \cup \sset 0$
is a $\Gal(F\unram/F)$-stable field satisfying
($\text{\textbf{CF}}_{F\unram}$)
that contains $\ff_F'$
and
is contained in $F\unram_0$\,.

\begin{rk}
It is easy to verify that the group $\mc F(F)$ does not
depend on the choice of $F'$.
However, for $p = \chr F$, it \emph{does} depend on the
choice of $\ff_F'$ (and, of course, of \ol F).
Since $\ff_F'$ may fail to be unique (see Theorem 10(a) of
\cite{cohen:complete-local-rings}), so may $\mc F(F)$;
but this ambiguity seems unavoidable.
\end{rk}

Note that, regardless of the values of $p$ and $\chr F$,
we have $\mc F(F) = \mc F(E)$ for any
discretely valued algebraic extension $E/F$.

\subsection{Absolute semisimplicity and topological
unipotence:  definitions and basic results}

\begin{dn}
\label{defn:top-F-ss-unip}
An element $\gamma \in G$ is
\emph{topologically $F$-unipotent} (in $G$)
if it belongs to $\bG(E)_{0+}$
for some finite extension $E/F$.
It is \emph{absolutely $F$-semisimple} (in $G$) if it is
semisimple and its character values belong to $\mc F(F)$.
If the image of $\gamma$ in \tG is absolutely
$F$-semisimple (respectively, topologically $F$-unipotent),
then we will say that $\gamma$ is \emph{absolutely
$F$-semisimple modulo \mo N{}} (respectively,
\emph{topologically $F$-unipotent modulo \mo N{}}).
\end{dn}

Note that an absolutely $F$-semisimple element need not
belong to $G_0$\,, and a topologically
$F$-unipotent element need not belong to $G_{0+}$ (but see
Proposition \ref{prop:tame-tunip}).
We will show later (see Corollary \ref{cor:abs-F-ss-tame})
that an absolutely $F$-semisimple element must be $F$-tame.

\begin{rk}
\label{rem:power-of-F-ss-unip}
Any power of an absolutely $F$-semisimple modulo \mo N
(respectively,
topologically $F$-unipotent modulo \mo N) element is again absolutely
$F$-semisimple modulo \mo N (respectively, topologically $F$-unipotent
modulo \mo N).
\end{rk}

\begin{rk}
\label{rem:top-F-ss-unip-bdd}
It is clear that a topologically $F$-unipotent element of
$G$ is bounded.
The character values of an absolutely $F$-semisimple element
lie in $\mc F(F) \subseteq (F\unram)\cross_0$\,;
so, by Remark \ref{rem:torus-filt}, $\gamma$ is
bounded.
It is an easy consequence that an element which is
absolutely $F$-semisimple or topologically $F$-unipotent
modulo \mo N is bounded modulo \mo N.
\end{rk}

\begin{lm}
\label{lem:bounded-and-top-F-unip}
If $\gamma, \delta \in G$ commute,
$\gamma$ is bounded, and
$\delta$ lies in $G_{0+}$\,, then there is a point
$x \in \BB(\bG, F)$ such that $\gamma\dota x = x$
and $\delta \in G_x^+$\,.
\end{lm}

\begin{proof}
Let $\varepsilon$ be a positive real number such that
$\mc S := \set{x \in \BB(\bG, F)}
	{\delta \in G_{x, \varepsilon}}$
is non-empty.
(Here, $G_{x, \varepsilon}$ is
a Moy--Prasad filtration subgroup.
See \cite{moy-prasad:k-types}*{\S 2.6} and
\cite{moy-prasad:jacquet}*{\S 3.2}.)
Certainly, \mc S is also closed and convex.
If $x \in \mc S$, then
$\delta = \lsup\gamma\delta
\in G_{\gamma\dota x, \varepsilon}$\,,
so $\gamma\dota x \in \mc S$.
That is, \mc S is $\gamma$-stable.
By Remark \ref{rem:equiv-bdd}, $\gamma$ fixes a point of
\mc S.
\end{proof}

\begin{lm}
\label{lem:top-F-unip-Jordan}
Suppose that $\gamma \in G$
has (ordinary) Jordan decomposition
$\gamma = \gamma\semi\gamma\unip$, and that
$\gamma\semi, \gamma\unip \in G$.
Then $\gamma$ is topologically $F$-unipotent modulo \mo N if
and only if $\gamma\semi$ is.
\end{lm}

\begin{proof}
Clearly, it suffices to prove this in case \mo N is the
trivial subgroup.

Suppose that
$\delta, \delta_+ \in G$ commute,
$\delta$ is topologically $F$-unipotent,
and $\delta_+$ is unipotent.
There is a finite extension $E/F$ such
that $\delta \in \bG(E)_{0+}$\,;
say $z \in \BB(\bG, E)$ is such that
$\delta \in \bG(E)_z^+$\,.
By Remark \ref{rem:equiv-unip}, there is a finite separable
extension $K/E$ such that $\delta_+ \in \bG(K)_{0+}$\,.
By Lemma \xref{exp-lem:field-descent} of
\cite{adler-spice:good-expansions},
$\delta \in \bG(K)_z^+$\,.
In particular, $\delta$ is bounded,
so,
by Lemma \ref{lem:bounded-and-top-F-unip}, there is a point
$y \in \BB(\bG, K)$ such that $\delta\dota y = y$ and
$\delta_+ \in \bG(E)_y^+$\,.
By Lemma \xref{exp-lem:stab-deep} of
\cite{adler-spice:good-expansions},
$\delta \in \bG(E)_y$\,.
By Lemma \xref{exp-lem:unipotent} of \emph{loc.\ cit.},
for $x \in (y, z)$ sufficiently close to $y$, we have that
$\delta \in \bG(E)_x^+$\,.
If, in addition, $x$ is so close to $y$ that it is contained
in a facet whose closure contains $y$, then
$\bG(E)_y^+ \subseteq \bG(E)_x^+$\,, so
$\delta\delta_+ \in \bG(E)_x^+$\,.
That is, $\delta\delta_+$ is topologically $F$-unipotent.

If we take $\delta = \gamma\semi$ and
$\delta_+ = \gamma\unip$, then we see that the topological
$F$-unipotence of $\gamma\semi$ implies that of $\gamma$.
If we take $\delta = \gamma$ and $\delta_+ = \gamma\unip\inv$,
then we see that the topological $F$-unipotence of $\gamma$
implies that of $\gamma\semi$.
\end{proof}

\begin{rk}
\label{rem:cent-of-ss}
Suppose that $\gamma \in G$ has (ordinary) Jordan
decomposition $\gamma = \gamma\semi\gamma\unip$.
Put $\wtilde\bH = C_\bG(\gamma\semi)$
and $\bH = \smash{\wtilde\bH}\conn$.
\begin{enumerate}
\item If $\chr F = 0$, then $\gamma\semi \in G$.
By Propositions 1.2(a), 9.1(1), and 13.19 of \cite{borel:linear},
\bH is a connected reductive $F$-group.
Certainly, $\gamma\semi \in H$ and
$\gamma\unip \in \wtilde H$.
Since the image of $\gamma\unip$ in the component group
$(\wtilde\bH/\bH)(F)$ is unipotent and has finite order, it
is trivial.  That is, $\gamma\unip \in H$,
so $\gamma \in H$.
\item If $\chr F > 0$, then, by
\cite{borel:linear}*{\S 4.1(a)}, there is some
$a \in \Z_{\ge 0}$ such that $\gamma\unip^{p^a} = 1$.
Then $\gamma\semi^{p^a} = \gamma^{p^a} \in G$.
Since an easy $\GL_n$ calculation shows that
$\wtilde\bH = C_\bG(\gamma\semi^{p^a})$,
hence that
$\bH = C_\bG(\gamma\semi^{p^a})\conn$,
we have again that \bH is a connected reductive $F$-group.
We have
$\gamma^{p^a} = \gamma\semi^{p^a} \in \bH(\ol F) \cap G = H$.
\end{enumerate}
\end{rk}

\begin{lm}
\label{lem:top-unip}
Suppose that $p > 0$.
Then an element of $G$ is topologically $F$-unipotent
modulo \mo N if and only if it is topologically
$p$-unipotent modulo $N$.
\end{lm}

\begin{proof}
Recall that, for every finite extension
$E/F$, Moy and Prasad have defined (in
\cite{moy-prasad:k-types}*{\S 2.6} and
\cite{moy-prasad:jacquet}*{\S 3.2}), for each
$x \in \BB(\bG, E)$,
an exhaustive filtration
$(\bG(E)_{x, r})_{r \in \R_{\ge 0}}$ of $\bG(E)$ by
subgroups such that $\bG(E)_{x, 0} = \bG(E)_x$
and $\bG(E)_{x, \varepsilon} = \bG(E)_x^+$
for sufficiently small positive $\varepsilon$.
Since $\bG(E)_{x, r}/\bG(E)_{x, r+}$ is a $p$-group for
$(x, r) \in \BB(\bG, E) \times \R_{> 0}$, we have that
a topologically $F$-unipotent modulo \mo N element is
topologically $p$-unipotent modulo $N$.

If $\gamma \in G$ is topologically $p$-unipotent modulo $N$,
then its image in \tG is topologically
$p$-unipotent.  Thus it suffices to prove that, if $\gamma$
is topologically $p$-unipotent, then it is topologically
$F$-unipotent.
Let $\gamma\semi$ and $\gamma\unip$ be the semisimple and
unipotent parts, respectively, of the (ordinary) Jordan
decomposition of $\gamma$; and put
$\bH = C_\bG(\gamma\semi)\conn$.

If $\chr F = 0$, then, by Remark \ref{rem:cent-of-ss},
we have $\gamma, \gamma\unip \in H$.
By Remark \ref{rem:equiv-unip},
there is a finite separable extension $E/F$ such that
\bH is $E$-split and
$\gamma\unip \in \bH(E)_{0+}$\,;
say $x \in \BB(\bH, E)$ is such that
$\gamma\unip \in \bH(E)_x^+$\,.

If $\chr F > 0$, then, again by Remark \ref{rem:cent-of-ss},
there is $a \in \Z_{\ge 0}$ such that
$\gamma^{p^a} \in H$
and $\gamma\unip^{p^a} = 1$.
In this case, let $E/F$ be a finite separable extension such
that \bH is $E$-split, and $x$ any point of $\BB(\bH, E)$.

In either case,
$\gamma^{p^n} \in \gamma\semi^{p^n}\bH(E)_x^+$
for all sufficiently large integers $n$.
Since also $\gamma^{p^n} \in \bH(E)_x^+$ for all
sufficiently large $n$ (by topological $p$-unipotence), we
have that $\gamma\semi^{p^n} \in \bH(E)_x^+$ for some
$n \in \Z_{\ge 0}$.
Let \bT be an $E$-split maximal torus in \bH (hence in \bG)
such that
$x$ belongs to the apartment of \bT in $\BB(\bH, E)$.
By Lemma \xref{exp-lem:levi-descent} of
\cite{adler-spice:good-expansions},
we have that $\gamma\semi^{p^n} \in \bT(E)_{0+}$\,.
Let $K/E$ be a finite extension such that
$\gamma\semi \in \bG(K)$.
By Remark \ref{rem:torus-filt}, the
character values of $\gamma\semi^{p^n}$
lie in $E\cross_{0+} \subseteq K\cross_{0+}$\,,
so the character valies of $\gamma\semi$ lie in
$K\cross_{0+}$\,,
so $\gamma\semi \in \bT(K)_{0+}$\,.
By another application of Lemma \xref{exp-lem:levi-descent}
of \emph{loc.\ cit.},
$\gamma\semi \in \bG(K)_{0+}$\,.
Thus $\gamma\semi$ is topologically $K$-unipotent.
By Lemma \ref{lem:top-F-unip-Jordan},
$\gamma$ is topologically $K$-unipotent.
It is then clear from the definition (see Definition
\ref{defn:top-F-ss-unip}) that it is
topologically $F$-unipotent.
\end{proof}

\begin{rk}
\label{rem:top-F-ss-unip-ascent}
Let $E/F$ be a discretely valued algebraic extension.
Since $\mc F(F) = \mc F(E)$,
an element of $G$ is absolutely $F$-semisimple modulo \mo N if
and only if it is absolutely $E$-semisimple modulo \mo N.
If $p = 0$, then Lemma \xref{exp-lem:domain-field-ascent} of
\cite{adler-spice:good-expansions} shows that an element of
$G$ is topologically $F$-unipotent modulo \mo N if and only
if it is topologically $E$-unipotent modulo \mo N.
If $p > 0$, then, by Lemma \ref{lem:top-unip}, the
topological $F$-unipotence and topological $E$-unipotence of
an element of $G$ are both equivalent to its topological
$p$-unipotence, hence to one another.
Thus an element of $G$ is
topologically $F$-unipotent modulo \mo N if and only if it
is topologically $E$-unipotent modulo \mo N.
\end{rk}

Our definition of topological $F$-Jordan decompositions is
almost the analogue one would expect of the definition of a
topological $p$-Jordan decomposition (see Definition
\ref{defn:top-p-Jordan}), except for one somewhat surprising
condition about tori.
Proposition \ref{prop:top-F-Jordan-defn} will show that this
condition can be omitted.

\begin{dn}
\label{defn:top-F-Jordan}
A \emph{topological $F$-Jordan decomposition modulo \mo N{}}
of an element $\gamma \in G$
is a pair of commuting elements $(\gamma\tsemi, \gamma\tunip)$
of $G$ such that
\begin{itemize}
\item the images of $\gamma\semi$ and $\gamma\tsemi$ in
$\tbG(\ol F)$ belong to a common \ol F-torus there,
\item $\gamma = \gamma\tsemi\gamma\tunip$,
\item $\gamma\tsemi$ is absolutely $F$-semisimple modulo \mo N,
and
\item $\gamma\tunip$ is topologically $F$-unipotent modulo \mo N.
\end{itemize}
We will sometimes just say that
$\gamma = \gamma\tsemi\gamma\tunip$ is a topological $F$-Jordan
decomposition modulo \mo N.
If \mo N is the trivial subgroup, then we will omit
``modulo \mo N''.
\end{dn}

\subsection{Uniqueness of topological Jordan decompositions}

\begin{pn}
\label{prop:unique-top-F-Jordan}
An element of $G$ has at most one topological
$F$-Jordan decomposition.
\end{pn}

\begin{proof}
Suppose that
$\gamma\tsemi\gamma\tunip = \gamma
= \gamma\tsemi'\gamma\tunip'$ are two topological $F$-Jordan
decompositions of an element $\gamma \in G$.
By Remark \ref{rem:top-F-ss-unip-ascent},
they remain topological $F$-Jordan decompositions if we
replace $F$ by a finite extension, so
we do so whenever necessary.

Since $\gamma\tsemi$ and $\gamma\tsemi'$
are semisimple and commute with $\gamma\tunip$ and
$\gamma\tunip'$, respectively, we have that
$(\gamma\tunip)\unip = \gamma\unip = (\gamma\tunip')\unip$.
Upon replacing $F$ by a finite extension, we may, and hence
do, assume that $\gamma\semi$
(hence also $(\gamma\tunip)\semi$ and $(\gamma\tunip')\semi$)
lie in $G$.
By Lemma \ref{lem:top-F-unip-Jordan},
$(\gamma\tunip)\semi$ and $(\gamma\tunip')\semi$ are
topologically $F$-unipotent.
Thus
$\gamma\tsemi(\gamma\tunip)\semi = \gamma\semi
= \gamma\tsemi'(\gamma\tunip')\semi$
are two topological $F$-Jordan decompositions of
$\gamma\semi$,
so we may, and hence do, assume that $\gamma$ is semisimple.

We show that $\gamma$, $\gamma\tsemi$, and $\gamma\tsemi'$
(hence also $\gamma\tunip$ and $\gamma\tunip'$) lie in a common
torus.
Let \bT and $\bT'$ be maximal \ol F-tori in \bG such that
$\gamma, \gamma\tsemi \in \bT(\ol F)$
and
$\gamma, \gamma\tsemi' \in \bT'(\ol F)$
(hence $\gamma\tunip' \in \bT'(\ol F)$).
Upon replacing $F$ by a finite extension,
we may, and hence
do, assume that $\bT'$ is $F$-split and
$\gamma\tunip' \in G_{0+}$\,.
By Lemma \xref{exp-lem:levi-descent} of
\cite{adler-spice:good-expansions}, we have
$\gamma\tunip' \in T'_{0+}$\,,
so that, by Remark \ref{rem:torus-filt},
the character values of $\gamma\tunip'$ lie in
$F\cross_{0+}$\,.
If $\alpha$ is a root of $\bT'$ in $C_\bG(\gamma)\conn$,
then $\alpha(\gamma) = 1$,
so
$\alpha(\gamma\tsemi') = \alpha(\gamma\tunip')\inv
\in \mc F(F) \cap F\cross_{0{+}}
= \sset 1$.
That is, $\gamma\tsemi'$ and $\gamma\tunip'$ are central in
$C_\bG(\gamma)\conn$,
hence belong to \bT.

Now the character values of $\gamma\tsemi$ and of
$\gamma\tsemi'$,
hence of $\gamma\tsemi'^{-1}\gamma\tsemi$, lie in $\mc F(F)$;
and those of $\gamma\tunip$ and of $\gamma\tunip'$,
hence of $\gamma\tunip'\gamma\tunip\inv$, lie in
$F\cross_{0{+}}$\,;
so, since $\mc F(F) \cap F\cross_{0+} = \sset 1$, we have that
$\gamma\tsemi'^{-1}\gamma\tsemi = \gamma\tunip'\gamma\tunip\inv$
equals $1$.
\end{proof}

It is an easy observation that, if $\bG'$ is a connected
reductive $F$-group and $f : \bG \to \bG'$ is an
$F$-morphism, then $f(\gamma)$ is absolutely $F$-semisimple
as long as $\gamma$ is.
Although we do not do so here, one can formulate a
condition on $f$ such that $f(\gamma)$ is topologically
$F$-unipotent as long as $\gamma$ is.
We record three consequences.

\begin{lm}
\label{lem:top-F-Jordan-central}
If $\gamma = \gamma\tsemi\gamma\tunip$ is a topological
$F$-Jordan decomposition, then
$\gamma\tsemi, \gamma\tunip \in Z(C_G(\gamma))$.
\end{lm}

\begin{proof}
Fix $g \in C_G(\gamma)$.
Then
$\gamma = \lsup g\gamma
= (\lsup g\gamma\tsemi)(\lsup g\gamma\tunip)$
is a topological $F$-Jordan decomposition.
By Proposition \ref{prop:unique-top-F-Jordan},
$\lsup g\gamma\tsemi = \gamma\tsemi$
and
$\lsup g\gamma\tunip = \gamma\tunip$.
\end{proof}

\begin{lm}
\label{lem:E-Jordan-is-F-ratl}
Suppose that
\begin{itemize}
\item $\gamma \in G$,
\item $E/F$ is a discretely valued separable extension,
and
\item $\gamma = \gamma\tsemi\gamma\tunip$ is a topological
$E$-Jordan decomposition.
\end{itemize}
Then $\gamma\tsemi, \gamma\tunip \in G$, and
$\gamma = \gamma\tsemi\gamma\tunip$ is a
topological $F$-Jordan decomposition.
\end{lm}

\begin{proof}
By Remark \ref{rem:top-F-ss-unip-ascent},
$\gamma = \gamma\tsemi\gamma\tunip$ is a topological
$\wtilde E$-Jordan decomposition, where $\wtilde E/F$ is the
Galois closure of $E/F$.
Then $\gamma = \sigma(\gamma\tsemi)\sigma(\gamma\tunip)$ is
also a topological $\wtilde E$-Jordan decomposition for
$\sigma \in \Gal(\wtilde E/F)$.
By Proposition \ref{prop:unique-top-F-Jordan},
$\gamma\tsemi, \gamma\tunip
\in \bG(\wtilde E)^{\Gal(\wtilde E/F)} = G$.
The last statement follows from another application of
Remark \ref{rem:top-F-ss-unip-ascent}.
\end{proof}

\begin{cor}
With the notation and hypotheses of Lemma
\ref{lem:E-Jordan-is-F-ratl}, suppose that $g \in \bG(E)$ is
such that $\lsup g\gamma \in G$.
Then $\lsup g\gamma\tsemi, \lsup g\gamma\tunip \in G$,
and
$\lsup g\gamma = (\lsup g\gamma\tsemi)(\lsup g\gamma\tunip)$
is a topological $F$-Jordan decomposition.
\end{cor}

\begin{proof}
It is clear that
$\lsup g\gamma = (\lsup g\gamma\tsemi)(\lsup g\gamma\tunip)$
is a topological $E$-Jordan decomposition.  Now the result
is an immediate consequence of Lemma
\ref{lem:E-Jordan-is-F-ratl}.
\end{proof}

\subsection{Relationship between algebraic and abstract
groups}

We now relate the abstract setting of \S\ref{sec:profinite}
to our present setting.
We have already seen that topological $F$- and
$p$-unipotence are equivalent when $p > 0$ (see Lemma
\ref{lem:top-unip}).  We prove below the analogous
results for absolute $F$-semisimplicity and topological
$F$-Jordan decompositions; but note that the formulations
are slightly more complicated.

\begin{lem}
\label{lem:abs-ss}
Suppose that $p > 0$ and $\gamma \in G$.
If $\gamma$ is absolutely $p$-semisimple modulo $N$, then it
is absolutely $F$-semisimple modulo \mo N.
If $\gamma$ is absolutely $F$-semisimple modulo \mo N,
then it is absolutely $p$-semisimple modulo $N$ if and only
if some finite power of it is topologically $p$-unipotent
modulo $N$.
\end{lem}

\begin{proof}
Suppose that $\gamma$ is absolutely $p$-semisimple modulo
$N$.
Certainly, some finite power of it is topologically
$p$-unipotent modulo $N$ (in fact, lies in $N$).
Further, its image \ogamma in \tG has finite, prime-to-$p$ order,
say $M$, hence is semisimple.
(Indeed, the unipotent part $(\ogamma)\unip$ of \ogamma also has
finite, prime-to-$p$ order.
If $\chr F = 0$ (respectively, $\chr F > 0$), then any non-trivial
unipotent element has infinite (respectively, $p$-power) order;
so $(\ogamma)\unip = 1$.)
By Lemma \xref{exp-lem:complete-subfield} of
\cite{adler-spice:good-expansions}, there is a complete
subfield $F'$ of $F$ such that
\begin{itemize}
\item $F/F'$ is unramified,
\item \tbG is defined over $F'$,
and
\item $\gamma \in \tbG(F')$.
\end{itemize}
Let
\begin{itemize}
\item $\wtilde\bT$ be a maximal $F'$-torus in \tbG such that
$\ogamma \in \wtilde\bT(F')$,
\item $E/F'$ the splitting field of $\wtilde\bT$,
and
\item $a$ an integer such that $a p \equiv 1 \pmod M$.
\end{itemize}
Then, for $\chi \in \bX^*(\bT)$,
we have
$\chi(\ogamma) = \chi(\ogamma^{a^m})^{p^m}
\in (E\cross)^{p^m}$
for all $m \in \Z_{\ge 0}$,
so
$\chi(\ogamma) \in (E\cross)^{p^\infty} \subseteq \mc F(F)$.
That is, \ogamma is absolutely $F$-semisimple, so
$\gamma$ is absolutely $F$-semisimple modulo \mo N.

Suppose that $\gamma$ is absolutely $F$-semisimple
modulo \mo N, and
$M \in \Z_{> 0}$ is such that $\gamma^M$ is topologically
$p$-unipotent modulo $N$.
Write $M = p^m M'$, with $m \in \Z_{\ge 0}$ and $M' \in \Z_{> 0}$
such that $M'$ is coprime to $p$.  Then $\gamma^{M'}$ is
also topologically $p$-unipotent modulo $N$, hence
(by Lemma \ref{lem:top-unip})
topologically $F$-unipotent modulo \mo N.
On the other hand,
by Remark \ref{rem:power-of-F-ss-unip}, $\gamma^{M'}$ is
absolutely $F$-semisimple modulo \mo N.
By Proposition \ref{prop:unique-top-F-Jordan}, the image of
$\gamma^{M'}$ in \tG is trivial, so
$\gamma^{M'} \in N$;
that is, $\gamma$ is absolutely $p$-semisimple modulo $N$.
\end{proof}

\begin{cor}
\label{cor:abs-ss-loc-cpct}
If $F$ is an algebraic extension of a locally compact field,
then an element of $G$
is absolutely $F$-semisimple modulo \mo N if and only if
it is absolutely $p$-semisimple modulo $N$.
\end{cor}

\begin{proof}
The `if' direction is clear.
For the `only if' direction, suppose that $\gamma \in G$ is
absolutely $F$-semisimple modulo \mo N.
By Remark \ref{rem:top-F-ss-unip-bdd},
$\gamma$ is bounded modulo \mo N;
so, by Remark \ref{rem:equiv-bdd},
it fixes some point $x \in \BB(\tbG, F)$.
By Lemma \xref{exp-lem:complete-subfield} of
\cite{adler-spice:good-expansions}, there is a locally
compact subfield $F'$ of $F$ such that
\begin{itemize}
\item \bG and \mo N are defined over $F'$,
\item $\gamma \in \bG(F')$,
and
\item $x \in \BB(\tbG, F')$.
\end{itemize}
Then the image of $\gamma$ in $\tbG(F')$ lies in
$\stab_{\tbG(F')}(x)$.  By Remarks \ref{rem:equiv-bdd} and
\ref{rem:bdd=cpct}, $\stab_{\tbG(F')}(x)$ is bounded, hence compact;
so its open subgroup $\tbG(F')_x^+$ has finite index.
That is, some (finite) power of $\gamma$ is
topologically $p$-unipotent modulo $N$.
Now the result follows from Lemma \ref{lem:abs-ss}.
\end{proof}

\begin{lm}
\label{lem:top-Jordan}
Suppose $p > 0$ and $\gamma\tsemi, \gamma\tunip \in G$.
If $(\gamma\tsemi, \gamma\tunip)$ is a topological
$p$-Jordan decomposition modulo $N$, then it is a
topological $F$-Jordan decomposition modulo \mo N.
If it is a topological $F$-Jordan decomposition modulo \mo N,
then it is a topological $p$-Jordan decomposition modulo $N$
if and only if $\gamma\tsemi$ is absolutely $p$-semisimple.
\end{lm}

\begin{proof}
It is clear from Lemma \ref{lem:top-unip} that,
if $(\gamma\tsemi, \gamma\tunip)$ is a topological
$F$-Jordan decomposition modulo \mo N, then it is a
topological $p$-Jordan decomposition modulo $N$ if and only
if $\gamma\tsemi$ is absolutely $p$-semisimple.

Suppose that $(\gamma\tsemi, \gamma\tunip)$ is a topological
$p$-Jordan decomposition modulo $N$.
Put $\gamma = \gamma\tsemi\gamma\tunip$.
By Lemmata \ref{lem:top-unip} and \ref{lem:abs-ss},
$\gamma\tsemi$ is absolutely $F$-semisimple modulo \mo N
and $\gamma\tunip$ is topologically $F$-unipotent \mo N.
Let $E/F$ be a finite extension such
that $\gamma\semi$ (hence $(\gamma\tunip)\semi$) belongs to
$\bG(E)$.
By Lemmata \ref{lem:top-F-unip-Jordan} and \ref{lem:top-unip},
$(\gamma\tunip)\semi$ is topologically $p$-unipotent modulo
$\mo N(E)$.
Since $\gamma\tsemi$ commutes with $\gamma\tunip$, it
commutes also with $(\gamma\tunip)\semi$, so
$\gamma\semi = \gamma\tsemi(\gamma\tunip)\semi$ is a
topological $p$-Jordan decomposition modulo $\mo N(E)$.
By Proposition \ref{prop:unique-top-p-Jordan_closure},
the image of $\gamma\tsemi$ in $\tbG(E)$ belongs to any
maximal \ol F-torus containing the image there of
$\gamma\semi$.
Thus $\gamma = \gamma\tsemi\gamma\tunip$ is a topological
$F$-Jordan decomposition modulo \mo N.
\end{proof}

\begin{cor}
\label{cor:top-Jordan-loc-cpct}
If $F$ is an algebraic extension of a locally compact field,
then a topological $F$-Jordan decomposition modulo \mo N is a
topological $p$-Jordan decomposition modulo $N$, and
conversely.
\end{cor}

\subsection{Stabilisers and parahorics}

It is a minor inconvenience in our arguments that the
stabiliser of a point in $\BB(\bG, F)$ may be strictly larger
than the associated parahoric subgroup.  The next
result shows that, under some circumstances, we may bring an
element of the stabiliser into the parahoric by passing to a
tame extension.

\begin{lm}
\label{lem:tame-para-ascent}
Suppose that
\begin{itemize}
\item $x \in \BB(\bG, F)$,
\item $g \in \stab_G(x)$,
and
\item $g^n \in G_x$ for some $n \in \Z_{> 0}$ indivisible
by $p$.
\end{itemize}
Then there exists a finite tame extension $L/F$ such that
$g \in \bG(L)_x$\,.
\end{lm}

\begin{proof}
Upon replacing $F$ by the splitting field of a maximal
$F$-tame torus, we may, and hence do, assume that \bG is
$F$-quasisplit.
Let \bS be a maximal $F$-split (hence maximal $F$-tame)
torus in \bG such that
$x$ belongs to the apartment $\AA(\bS, F)$ of \bS,
and \bT the maximal $F$-torus in \bG containing \bS.
If $a$ is a root of \bS in \bG, then let $\alpha$ be a root
of \bT in \bG restricting to $\alpha$, and write
$F_a$ for the fixed field in $F\sep$ of
$\stab_{\Gal(F\sep/F)}(\alpha)$
(the field denoted by $L_a$ in
\cite{bruhat-tits:reductive-groups-2}*{\S\S 4.1.8 and 4.1.14}).
Up to $F$-isomorphism, this field does not depend on the
choice of $\alpha$.
By replacing $F$ by a further finite tame extension if necessary,
we may, and hence do, assume that
$F$ contains the $n$th roots of unity in $F\sep$, and
all the extensions $F_a/F$ are totally wildly ramified.
Let $L'/F$ be a totally (tamely) ramified extension of
degree $n$.
Note that \bS is still a maximal $L'$-split torus.
Fix a root $a$ of \bS in \bG.
With the obvious notation, $L'_a$ is the fixed field in
$F\sep$ of
$\stab_{\Gal(F\sep/L')}(\alpha)
= \Gal(F\sep/L') \cap \stab_{\Gal(F\sep/F)}(\alpha)$
---
that is, $L'_a = L' F_a$.
Since $F_a/F$ is totally wildly ramified
and $L'/F$ is tamely ramified,
it follows that
$L'_a/F_a$ is a totally ramified extension of degree $n$.

Choose a chamber $C$ in $\AA(\bS, F)$ containing $x$ in its
closure,
and a special vertex $o$ in the closure of $C$.
By regarding $o$ as an origin, we may, and hence do,
identify $\AA(\bS, F)$ with $\bX_*(\bS) \otimes_\Z \R$,
hence the affine $F$-roots on $\AA(\bS, F)$
(in the sense of \cite{moy-prasad:k-types}*{\S 2.5},
not \cite{adler-spice:good-expansions}*{\S 2.2}) with
certain functions on $\bX_*(\bS) \otimes_\Z \R$ of the form
$y \mapsto \langle a, y\rangle + r$
with $a$ a root of \bS in \bG and $r \in \R$.
(Here, $\langle\cdot, \cdot\rangle$ is the usual pairing
between $\bX^*(\bS)$ and $\bX_*(\bS)$.)
Specifically, $r$ must belong to the set denoted by
$\Gamma'_a$ in
\cite{bruhat-tits:reductive-groups-1}*{\S 6.2.2}.
By \cite{bruhat-tits:reductive-groups-2}*{\S 4.2.21}
(adapted to our choice of
origin, which is different from the one in \S 4.2.2 of
\emph{loc.\ cit.}),
we have $\Gamma'_a = \ord(F_a\cross)$.
Denote by $\lsub F\mc H$ the collection of zero-sets of
affine $F$-roots.

We have that
$\lsub F W\textsub{aff} := N_G(T)/T\subb$\,, viewed as a group of
affine transformations of $\AA(\bS, F)$, is isomorphic to the
semi-direct product $\lsub F\Lambda \rtimes \lsub F W$, where
$\lsub F\Lambda = T/T\subb$
is a lattice of translations, and
$\lsub F W \cong N_G(T)/T$ is
the finite group generated by the reflections
through the hyperplanes in $\lsub F\mc H$
passing through $o$.
Let $\lsub F W'$ be the (normal) subgroup of $\lsub F W\textsub{aff}$
generated by reflections through the hyperplanes in $\lsub F\mc H$.
Then $\lsub F W' \cap \lsub F\Lambda$ is generated by
translations by elements of the form
$\gamma' a^\vee$,
where $a$ is a root of \bS in \bG,
$a^\vee$ is the associated coroot,
and $\gamma' \in \ord(F_a\cross)$.
We will denote by a left subscript $L'$ the analogues
over $L'$ of the objects defined over $F$ above.
Then the fact that
$\ord({L'_a}\cross) = \tfrac 1 n\ord(F_a\cross)$
and the
obvious analogue for $L'$ of our discussion above
for $F$ show that
\begin{equation}
\tag{$*$}
\text{
if $\tau \in \lsub F\Lambda$ satisfies
$\tau^n \in \lsub F W'$, then
$\tau \in \lsub{L'}W'$\,.
}
\end{equation}

Let $\Omega$ be the image of $\sset x$ in the reduced
building $\rBB(\bG, F)$,
and $f = f_\Omega'$ the optimisation of the function $f_\Omega$ of
\cite{bruhat-tits:reductive-groups-2}*{\S 4.6.26}.
Then,
by Proposition 4.6.28(i) and D\'efinition 5.2.6 of
\emph{loc.\ cit.},
the group of integer points of the scheme $\mf G_f^0$ of
\S 4.6.2 of \emph{loc.\ cit.}\ is the parahoric $G_x$\,.
By Corollaire 4.6.12 of \emph{loc.\ cit.},
there exists, for each root $a \in \Phi_f$, an affine
transformation $w_a \in W_f$ such that the linear
part of $w_a$ is reflection in the zero-set of $a$;
and $W_f$ is generated by the elements $w_a$.
Here, $\Phi_f$ is the set of gradients of affine $F$-roots
vanishing at $x$, and $W_f$ is as in
4.6.3(6) of \emph{loc.\ cit.}
Fix $a \in \Phi_f$, and let $\psi$ be the affine $F$-root
with gradient $a$ that vanishes at $x$.
Since $w_a$ fixes $x$, it must actually be reflection
in the zero-set of $\psi$.
That is, $W_f$ is generated by the reflections
through hyperplanes in $\lsub F\mc H$ passing through $x$.
By Proposition V.3.2 of \cite{bourbaki:lie-gp+lie-alg_4-6},
it is actually the stabiliser of $x$ in $\lsub F W'$.
Since $N_G(T) \cap G_x = N_f^0$, in the notation of 4.6.3(5)
of \emph{loc.\ cit.}, and since $W_f$ is the image in
$\lsub F W\textsub{aff}$ of $N_f^0$, we have that
\begin{equation}
\tag{$**_F$}
\text{
the image of $N_G(T) \cap G_x$ in
$\lsub F W\textsub{aff}$ is the stabiliser in
$\lsub F W'$ of $x$.
}
\end{equation}
Of course, there is an analogous statement, which we will
denote by ($**_{L'}$), when $F$ is replaced by $L'$.

By Proposition 4.6.28 of
\cite{bruhat-tits:reductive-groups-2},
$N_G(T) \cap g G_x \ne \emptyset$.
We may, and hence do,
replace $g$ by an element of this intersection.
Then write $w(g)$ for the image of $g$ in $\lsub F W\textsub{aff}$,
and let $\tau \in \lsub F\Lambda$ be such that
$w(g) \in \tau\dotm\lsub F W'$.
Then $w(g)^n \in \tau^n\dotm\lsub F W'$.
Since $g^n \in N_G(T) \cap G_x$\,,
we have by ($**_F$) that $w(g)^n \in \lsub F W'$.
Thus, $\tau^n \in \lsub F W'$.
By ($*$), we have that
$\tau \in \lsub{L'}W'$, so
\begin{equation}
\tag{$*{*}*$}
w(g) \in \tau\dotm\lsub F W' \subseteq \lsub{L'}W'.
\end{equation}

Since $w(g)$ stabilises $x$, we have by
($**_{L'}$) and ($*{*}*$) that it belongs to
the image of $N_\bG(\bT)(L') \cap \bG(L')_x$
in $\lsub{L'}W\textsub{aff}$\,.
That is,
$g \in (N_\bG(\bT)(L') \cap \bG(L')_x)\bT(L')\subb$\,.
In particular, $\bG(L')_x\dotm g$ contains an element
of $\bT(L')$, say $t$.
Then $t^n \in \bG(L')_x$\,.
By Lemma \xref{exp-lem:levi-descent} of
\cite{adler-spice:good-expansions}, we have that 
$t^n \in \bT(L')_0$\,.
Now we imitate the proof of Lemma
\xref{exp-lem:torus-field-descent} of \emph{loc.\ cit.}\
to show that there is a
finite tame extension $L/L'$ such that $t \in \bT(L)_0$\,.
Denote by $M$ a totally ramified extension of
${L'}\unram$ of degree $n$.
Then, in the notation of \cite{kottwitz:isocrystals-2}*{\S 7.3}
(except that our $M$ and $L'$ are Kottwitz's $L'$ and $L$,
respectively; so $\beta$ is the inclusion of $\bT(L')$ in
$\bT(M)$), we have by (7.3.2) of \emph{loc.\ cit.}\ that
$$
\alpha(w_{\bT(M)}(\beta(t))) = \alpha(N(w_{\bT(L')}(t)))
= n w_{\bT(L')}(t) = w_{\bT(L')}(t^n).
$$
By Lemma 2.3 of \cite{rapoport:T1-is-T0},
$\bT(L')_0 = \ker w_{\bT(L')}$
and
$\bT(M)_0 = \ker w_{\bT(M)}$\,.
In particular, $\alpha(w_{\bT(M)}(\beta(t))) = 0$, so, since
$\alpha$ is an injection, $t = \beta(t) \in \bT(M)_0$\,.
Now let $L/L'$ be any finite subextension of $M/L'$ such
that $M/L$ is unramified.
Then $t \in \bT(M)_0^{\Gal(M/L)} = \bT(L)_0$\,.

By Lemma \xref{exp-lem:levi-descent} of
\cite{adler-spice:good-expansions}, we have that
$t \in \bG(L)_x$\,,
so
$g \in \bG(L)_x\dotm t = \bG(L)_x$\,.
\end{proof}

\subsection{Existence of topological Jordan decompositions}
\label{subsec:existence}

The following two results show that the answers are ``yes''
to the analogues of the questions posed in
\cite{moy:displacement}*{\S\S 5.7 and 5.10},
where semisimplicity and unipotence are replaced by absolute
$F$-semisimplicity and topological $F$-unipotence.
We must impose at first a somewhat artifical tameness hypothesis,
but Corollary \ref{cor:abs-F-ss-tame} below will show that
it can be omitted.

\begin{pn}
\label{prop:x-depth}
If $\gamma$ is absolutely semisimple and $F$-tame, then
$\BB(C_\bG(\gamma), F)
= \set{x \in \BB(\bG, F)}{\gamma\dota x = x}$.
\end{pn}

\begin{proof}
Denote the right-hand set above by $\BB(\gamma)$.
By Proposition \xref{exp-prop:compatibly-filtered-tame-rank}
of \cite{adler-spice:good-expansions}, we have that
$C_\bG(\gamma)$ is a compatibly filtered $F$-subgroup of
\bG, in the sense of Definition
\xref{exp-defn:compatibly-filtered} of \emph{loc.\ cit.}
In particular, $\BB(C_\bG(\gamma), E)$ may be regarded
canonically as a subset of $\BB(\bG, E)$
for all discretely valued tame extensions $E/F$,
so that the statement makes sense.
Since
$\BB(C_\bG(\gamma), F)
= \BB(C_\bG(\gamma), E)^{\Gal(E/F)}$
for any discretely valued,
tame, Galois extension $E/F$, we may, and hence do,
assume that $F$ is strictly Henselian (hence that \bG is
$F$-quasisplit)
and that $\gamma$ belongs to a maximal $F$-split torus 
\bS in \bG.
Since $\gamma$ is bounded (by Remark
\ref{rem:top-F-ss-unip-bdd})
and $S\subb = S_0$ (by Remark \ref{rem:torus-filt}),
we have
$\gamma \in S_0 \subseteq G_0$
and $\BB(C_\bG(\gamma), F) \subseteq \BB(\gamma)$.

Suppose that $x \in \BB(C_\bG(\gamma), F)$, and
$y \in \BB(\gamma)$ lies in a facet of $\BB(\bG, F)$ whose
closure contains $x$.
Denote by $g \mapsto \ol g$ the reduction map
$G_x \to \ms G_x\conn(\ff)$.

We have that
$G_x^+ \subseteq G_y^+ \subseteq G_y \subseteq G_x$\,,
and the images in $\ms G_x\conn(\ff)$ of $G_y$ and $G_y^+$
are the groups of \ff-points of a parabolic \ff-subgroup $\ms P_y$
and of its unipotent radical $\ms U_y$, respectively.
Let \ms T be the \ff-split maximal torus in $\ms G_x\conn$
such that the image of $S_0$ in $\ms G_x\conn(\ff)$ is
$\ms T(\ff)$.
By Lemma \xref{exp-lem:stab-deep} of
\cite{adler-spice:good-expansions},
we have that $\gamma \in G_y \subseteq G_x$\,,
so $\ogamma \in \ms T(\ff) \cap \ms P_y(\ff)$ is
semisimple.
Thus it lies in a maximal \ff-torus $\ms T'$
of $\ms P_y$.

We claim that there is an $F$-split torus $\bS'$ in
$C_\bG(\gamma)\conn$ such that $x$ lies in the apartment of
$\bS'$
and
the image of $S'_0$ in $\ms G_x\conn(\ff)$ is
$\ms T'(\ff)$.
Indeed,
since \bS is $F$-split, there is an isomorphism
$i : \bX^*(\bS) \to \bX^*(\ms T)$ such that, for all
$\chi \in \bX^*(\bS)$, the image in $\ff\cross$ of
$\chi(\gamma) \in F\cross_0$ is $i(\chi)(\ogamma)$.
By Proposition 3.5.4 of \cite{steinberg:conjugacy},
$C_{\ms G_x}(\ogamma)\conn(\ff)$ is generated by
$\ms T(\ff)$ and the \ff-points of those root subgroups
corresponding to roots of \ms T in $\ms G_x\conn$ that
vanish at \ogamma.
Let $\ol\alpha$ be such a root.
By Corollaire 4.6.12(i) of
\cite{bruhat-tits:reductive-groups-2} (applied to the
function $f = f_\Omega'$ occurring in the proof of Lemma
\ref{lem:tame-para-ascent}),
$\alpha := i\inv(\ol\alpha)$ is a root of \bS in \bG.
Let $U \subseteq U_\alpha \cap G_x$ be the affine root subgroup of $G$
that maps onto the \ff-points of the root subgroup
$\ms U_{\ol\alpha}$ of $\ms G_x\conn$.
Since the image of $\alpha(\gamma)$ in $\ff\cross$ is
$\ol\alpha(\ogamma) = 1$, we have that
$\alpha(\gamma) \in \mc F(F) \cap F\cross_{0+} = \sset 1$.
By Proposition 3.5.4 of \cite{steinberg:conjugacy},
$\alpha$ is a root of \bS in $C_\bG(\gamma)\conn$,
so $U \subseteq C_G(\gamma)\conn \cap G_x$.
That is, the image in $\ms G_x\conn(\ff)$ of
$C_G(\gamma)\conn \cap G_x$ includes
$C_{\ms G_x}(\ogamma)\conn(\ff)$.
(Although we do not need to do so here, one can show that
the image is precisely $C_{\ms G_x\conn}(\ogamma)(\ff)$.)
Since \ms T and $\ms T'$ are maximal \ff-tori in
$C_{\ms G_x}(\ogamma)\conn$ and \ff is algebraically
closed, there is an element
$\ol c \in C_{\ms G_x}(\ogamma)\conn(\ff)$
such that $\ms T' = \lsup{\ol c}\ms T$.
Let $c \in C_G(\gamma)\conn \cap G_x$ be an element whose
image in $\ms G_x\conn(\ff)$ is \ol c.
Then $\bS' := \lsup c\bS$ certainly contains $x$ in its
apartment, and has the property that the image of
$S'_0$ in $\ms G_x\conn(\ff)$ is $\ms T'(\ff)$.

Note that we may, and hence do,
also regard $\ms T'$ as a torus in
$\ms P_y/\ms U_y = \ms G_y\conn$.
By Proposition 5.1.10 of
\cite{bruhat-tits:reductive-groups-2},
there is an $F$-split torus $\bS''$ in \bG
such that the apartment of $\bS''$ contains $y$
and the image of $S''_0$ in $\ms G_y\conn(\ff)$
is $\ms T'(\ff)$.
Since $y$ lies in a facet whose closure contains $x$,
the apartment of $\bS''$ also contains $x$.
By Proposition 4.6.28(iii) of
\cite{bruhat-tits:reductive-groups-2},
there is an element $k \in G_x$ such that
$\bS'' = \lsup k\bS'$.
Since $S''_0$ and $S'_0$ have the same image, namely
$\ms T'(\ff)$,
in $\ms G_x\conn(\ff)$, we have that
$\ol k \in N_{\ms G_x\conn}(\ms T')(\ff)$.
As in the proof of Lemma \ref{lem:tame-para-ascent}, one
sees from Corollaire 4.6.12(ii) of \emph{loc.\ cit.}\
that \ol k lies in the image in
$\ms G_x\conn(\ff)$ of $N_G(S') \cap G_x$\,.
Thus, there are
$k_+ \in G_x^+ \subseteq G_y^+$
and
$n \in N_G(S') \cap G_x$
such that $k = k_+ n$.
Then
$y = k_+\inv y
\in \AA(\lsup{k_+\inv}\bS'', F)
= \AA(\lsup n\bS', F) = \AA(\bS', F)
\subseteq \BB(C_\bG(\gamma), F)$.

We have shown that $\BB(C_\bG(\gamma), F)$ is open in
$\BB(\gamma)$.  Since it is a union of apartments, it is also
closed there.  Since
$\BB(\gamma)$ is connected (even convex), and since
$\BB(C_\bG(\gamma), F)$ is non-empty, we have the desired
equality.
\end{proof}

\begin{lm}
\label{lem:top-Jordan-to-Jordan}
Suppose that
\begin{itemize}
\item $x \in \BB(\bG, F)$,
\item $\gamma\tsemi, \gamma\tunip \in \stab_G(x)$
are absolutely $F$-semisimple and topologically
$F$-unipotent, respectively,
and
\item $\gamma\tsemi$ is $F$-tame.
\end{itemize}
Then the images of $\gamma\tsemi$ and $\gamma\tunip$ in
$\ms G_x(\ff)$ are semisimple and unipotent, respectively.
\end{lm}

\begin{proof}
We first show that the image of $\gamma\tunip$ is unipotent.
If $p > 0$, then we have by Lemma \ref{lem:top-unip}
that $\gamma\tunip$ is topologically $p$-unipotent, hence
that the image of $\gamma\tunip$ in $\ms G_x(\ff)$ has
$p$-power order.  By \cite{borel:linear}*{\S 4.1(a)}, it
is unipotent.

If $p = 0$, then let $E/F$ be a finite
extension such that $\gamma\tunip \in \bG(E)_{0+}$\,.
Since $E/F$ is tame, we have by Lemma
\xref{exp-lem:domain-field-ascent} of
\cite{adler-spice:good-expansions} that
$\gamma\tunip \in G_{0+}$\,.
Choose $z \in \BB(\bG, F)$ such
that $\gamma\tunip \in G_z^+$\,.
By Lemma \xref{exp-lem:stab-deep} of \emph{loc.\ cit.},
$\gamma\tunip \in G_x$\,.
By Lemma \xref{exp-lem:unipotent} of \emph{loc.\ cit.},
there is a point $y \in (x, z)$ such that
$\gamma\tunip \in G_y^+$ and $y$ belongs to a facet of
$\BB(\bG, F)$ whose closure contains $x$.
Then the image of $G_y^+$ in $G_x$ is the group of
\ff-points of the unipotent radical of a parabolic
\ff-subgroup of $\ms G_x\conn$.
In particular, the image of $\gamma\tunip$ is unipotent.

Now we show that the image of $\gamma\tsemi$ is semisimple.
Let $L/F\unram$ be a finite tame extension such that
$\gamma\tsemi$ belongs to an $L$-split torus.
By Lemma \xref{exp-lem:field-descent} of
\cite{adler-spice:good-expansions},
$\bG(L)_x^+ \cap \bG(F\unram) = \bG(F\unram)_x^+$\,,
so
$$
\ms G_x(\ol\ff) = \stab_{\bG(F\unram)}(x)/\bG(F\unram)_x^+
\subseteq \stab_{\bG(L)}(x)/\bG(L)_x^+ = \ms G_x^L(\ol\ff);
$$
that is, $\ms G_x$ is an \ol\ff-subgroup of
$\ms G_x^L$.
Thus we may, and hence do, assume, upon
replacing $F$ by $L$, that $\gamma\tsemi$ belongs to an
$F$-split torus.
By Proposition \ref{prop:x-depth}, we have that
there is a maximal $F$-split torus \bS
whose apartment contains $x$
such that $\gamma\tsemi \in S$.
Since $\gamma\tsemi$ is bounded (by Remark
\ref{rem:top-F-ss-unip-bdd})
and $S\subb = S_0$ (by Remark \ref{rem:torus-filt}),
we have $\gamma\tsemi \in S_0$\,.
Then the image of $\gamma\tsemi$ in $\ms G_x(\ff)$ belongs
to the image of $S_0$ there, which is the group of
\ff-rational points of an \ff-torus.
\end{proof}

\begin{lm}
\label{lem:fixed-pts}
If $\gamma = \gamma\tsemi\gamma\tunip$ is a topological
$F$-Jordan decomposition, then a point $x$ of $\BB(\bG, F)$ is
fixed by $\gamma$ if and only if it is fixed by
$\gamma\tsemi$ and $\gamma\tunip$.
\end{lm}

\begin{proof}
The `if' direction is obvious, so we need only prove the
`only if' direction.
By Remark \ref{rem:top-F-ss-unip-ascent}, it suffices to
prove this result over any finite
extension of $F$; so we may, and hence do, assume that
$\gamma\tsemi$ belongs to an $F$-split maximal torus in \bG.

Denote by $\BB(\gamma)$ the fixed points of $\gamma$, and
similarly for $\gamma\tsemi$ and $\gamma\tunip$.
Suppose that
$x \in \BB(\gamma\tsemi) \cap \BB(\gamma\tunip)
	\subseteq \BB(\gamma)$,
and $y \in \BB(\gamma)$ belongs to a facet whose closure
contains $x$.
Denote by $g \mapsto \ol g$ the reduction map
$G_x \to \ms G_x\conn(\ff)$.

The image of $G_y$ in $\ms G_x\conn(\ff)$ is the group of
\ff-points of a parabolic \ff-subgroup $\ms P_y$ of
$\ms G_x\conn$.
Since $\gamma$ normalizes $G_y$\,,
\ogamma normalizes $\ms P_y(\ff)$; so, by
Theorem 11.16 of \cite{borel:linear},
$\ogamma \in \ms P_y(\ff)$.
Then also
$(\ogamma)\semi \in \ms P_y(\ff)$
and
$(\ogamma)\unip \in \ms P_y(\ff)$.
By Lemma \ref{lem:top-Jordan-to-Jordan},
$(\ogamma)\semi = \ol{\gamma\tsemi}$
and
$(\ogamma)\unip = \ol{\gamma\tunip}$.
Since the preimage of $\ms P_y(\ff)$ in $G_x$ is $G_y$\,,
we have that $\gamma\tsemi$ and $\gamma\tunip$ lie in
$G_y$\,.
In particular,
$y \in \BB(\gamma\tsemi) \cap \BB(\gamma\tunip)$.
That is, $\BB(\gamma\tsemi) \cap \BB(\gamma\tunip)$ is open
in $\BB(\gamma)$.
Since it is also closed, and since $\BB(\gamma)$ is
connected (even convex), we have equality, as desired.
\end{proof}

Now we are in a position to prove an existence result
for topological $F$-Jordan decompositions
analogous to Proposition \ref{prop:intrinsic-tu}.
A more refined version of this result appears as
Theorem \ref{thm:top-F-Jordan-projects} below.

\begin{pn}
\label{prop:tJd=cpct}
An element $\gamma \in G$ has a topological $F$-Jordan
decomposition
$\gamma = \gamma\tsemi\gamma\tunip$
if and only if it is bounded.
In this case, $\gamma\tsemi$ is $F$-tame.
\end{pn}

\begin{proof}
Suppose that $\gamma = \gamma\tsemi\gamma\tunip$ is a
topological $F$-Jordan decomposition.
By Remark \ref{rem:top-F-ss-unip-bdd},
$\gamma\tsemi$ and $\gamma\tunip$ are bounded.
By Lemma \ref{lem:bounded-and-top-F-unip},
there is a point $x \in \BB(\bG, F)$ fixed by both,
so $\gamma\dota x = x$.
By Remark \ref{rem:equiv-bdd}, $\gamma$
is bounded.

Now suppose that $\gamma$ is bounded.
By Remark \ref{rem:top-F-ss-unip-ascent} and Lemma
\ref{lem:E-Jordan-is-F-ratl}, we may, and hence
do, replace $F$ by discretely valued tame extensions as
necessary.
In particular, we will assume throughout that $F = F\unram$.
Put $\bH = C_\bG(\gamma\semi)\conn$.

If $p > 0$, then, as in Remark \ref{rem:cent-of-ss},
let $a \in \Z_{\ge 0}$ be so large that
$\gamma^{p^a}$ and $\gamma\unip^{p^a}$
belong to $H$.
Let $E/F$ be a finite separable extension such that \bH is
$E$-split.
By Proposition \xref{exp-prop:compatibly-filtered-tame-rank}
of \cite{adler-spice:good-expansions},
\bH is a compatibly filtered $E$-subgroup of \bG,
in the sense of Definition
\xref{exp-defn:compatibly-filtered} of \emph{loc.\ cit.}
In particular, the building of $\BB(\bH, E)$ may be embedded
isometrically and $\gamma^{p^a}$-equivariantly into
$\BB(\bG, E)$.
Thus the orbits of $\gamma^{p^a}$ in $\BB(\bH, E)$,
hence in $\BB(\bH, F)$, are bounded;
that is, $\gamma^{p^a}$ is bounded (in $H$).
By Remark \ref{rem:equiv-bdd}, there is a point
$x \in \BB(\bH, F)$ fixed by $\gamma^{p^a}$.
Denote by $h \mapsto \ol h$ the reduction map
$\stab_H(x) \to \ms H_x(\ff)$.
Let $b \in \Z_{\ge 0}$ be so large that the order of
$\gamma^{p^{a + b}} \in \stab_H(x)$ modulo $H_x$ is
indivisible by $p$.
By Lemma \ref{lem:tame-para-ascent}, we may, and hence do,
assume, upon replacing $F$ by a finite tame extension, that
$\gamma^{p^{a + b}} \in H_x$\,.
By Remark \ref{rem:equiv-unip} and Lemma \ref{lem:top-unip},
$\gamma\unip^{p^{a + b}}$ is topologically $p$-unipotent.
Let $c \in \Z_{\ge 0}$ be so large that
$\gamma\unip^{p^{a + b + c}} \in H_x^+$\,,
and put $n = a + b + c$.
Since $\gamma\semi^{p^n} \in Z(H)$ and \ff is algebraically
closed, we have that
$\ol{\gamma\semi^{p^n}} \in Z(\ms H_x\conn(\ff))
	= Z(\ms H_x\conn)(\ff)$.
Let \ms T be a maximal \ff-torus (necessarily \ff-split)
in $\ms H_x\conn$,
so that $\ol{\gamma\semi^{p^n}} \in \ms T(\ff)$.
By Proposition 5.1.10 of
\cite{bruhat-tits:reductive-groups-2}, there exists a
maximal $F$-split torus \bS in \bH such that
$x \in \AA(\bS, F)$ and the image of $S_0$ in
$\ms H_x\conn(\ff)$ is $\ms T(\ff)$.
Let $\delta$ be a preimage in $S_0$ of
$\ol{\gamma\semi^{p^n}}$\,,
so that $\delta\inv\gamma\semi^{p^n} \in H_x^+$\,.
For $\chi \in \bX^*(\bS)$, let
$s_\chi'$ be the unique element of $\mc F(F)$ such that
$\chi(\delta) \equiv s_\chi' \pmod{F\cross_{0+}}$,
and $s_\chi$ the unique element of $\mc F(F)$ such that
$s_\chi^{p^n} = s_\chi'$.
Finally, let $\gamma\tsemi$ be the unique element of
$S$ such that $\chi(\gamma\tsemi) = s_\chi$ for all
$\chi \in \bX^*(\bS)$.
Clearly, $\gamma\tsemi$ is absolutely $F$-semisimple and
$F$-tame (even $F$-split, in the obvious language).
Put $(\gamma\semi)\tunip := \gamma\tsemi\inv\gamma\semi$.
By Remark \ref{rem:torus-filt},
$\gamma\tsemi^{-p^n}\delta
\in S_{0+} \subseteq H_x^+$\,.
Thus
$(\gamma\semi)\tunip^{p^n}
= \gamma\tsemi^{-p^n}\gamma\semi^{p^n} \in H_x^+$\,,
so $(\gamma\semi)\tunip$
is topologically $p$-unipotent.
By Lemma \ref{lem:top-unip},
$(\gamma\semi)\tunip$ is topologically $K$-unipotent
(where $K/F$ is a finite extension such that
$\gamma\semi \in \bH(K)$).
Thus,
$\gamma\semi = \gamma\tsemi(\gamma\semi)\tunip$
is a topological $K$-Jordan decomposition.
By Lemma \ref{lem:top-F-Jordan-central},
$\gamma\tsemi$ and $(\gamma\semi)\tunip$ commute with
$C_{\bG(K)}(\gamma\semi)$; in particular, with
$\gamma\unip$.
Put
$\gamma\tunip := \gamma\tsemi\inv\gamma
	= (\gamma\semi)\tunip\gamma\unip
\in G$.
Since
$\gamma\tunip^{p^n}
= (\gamma\semi)\tunip^{p^n}\gamma\unip^{p^n}
\in H_x^+$\,,
we have that $\gamma\tunip$ is topologically $p$-unipotent,
hence, by another application of Lemma \ref{lem:top-unip},
topologically $F$-unipotent.
Thus,
$\gamma = \gamma\tsemi\gamma\tunip$ is the desired
topological $F$-Jordan decomposition.

If $p = 0$, then $\gamma\unip \in H$.
By Remark \ref{rem:equiv-unip}, we may, and hence do,
assume, upon replacing $F$ by a finite (necessarily tame) extension,
that $\gamma\unip \in H_{0+}$ and \bH is $F$-split.
By Lemma \ref{lem:bounded-and-top-F-unip},
there is a point $x \in \BB(\bH, F)$ such that
$\gamma\dota x = x$ and $\gamma\unip \in H_x^+$\,.
Let \bT be an $F$-split maximal torus in \bH whose apartment
contains $x$.
Then $\gamma\semi \in T$ fixes $x$, hence is bounded.
By Remark \ref{rem:torus-filt},
the character values of $\gamma\semi$ lie in $F\cross_0$\,.
For $\chi \in \bX^*(\bT)$, let $s_\chi$ be the unique
element of $\mc F(F)$ such that
$\chi(\gamma\semi) \equiv s_\chi \pmod{F\cross_{0+}}$.
In particular, $s_\alpha = 1$ for all
roots $\alpha$ of \bT in \bH.
Let $\gamma\tsemi$ be the unique element of $T$ such
that $\chi(\gamma\tsemi) = s_\chi$ for all
$\chi \in \bX^*(\bT)$.
In particular, $\alpha(\gamma\tsemi) = 1$ for all
roots $\alpha$ of \bT in \bH, so
$\gamma\tsemi \in Z(H)$.
Clearly,
$\gamma\tsemi$ is $F$-tame and absolutely $F$-semisimple,
and belongs to an \ol F-torus containing $\gamma\semi$.
Moreover, by Remark \ref{rem:torus-filt},
$\gamma\tsemi\inv\gamma\semi
\in T_0^+ \subseteq H_x^+$\,.
Thus
$\gamma\tunip := \gamma\tsemi\inv\gamma
= (\gamma\tsemi\inv\gamma\semi)\gamma\unip \in H_x^+$\,.
By Proposition
\xref{exp-prop:compatibly-filtered-tame-rank}
of \cite{adler-spice:good-expansions},
\bH is a compatibly filtered $F$-subgroup of \bG,
in the sense of Definition
\xref{exp-defn:compatibly-filtered} of \emph{loc.\ cit.}
In particular, we may regard $x$ (non-canonically) as a
point of $\BB(\bG, F)$.
Then $H_x^+ \subseteq G_x^+$\,,
so $\gamma\tunip$ is topologically $F$-unipotent (in $G$).
Clearly, $\gamma\tunip \in H$ commutes with
$\gamma\tsemi \in Z(H)$.
Thus,
$\gamma = \gamma\tsemi\gamma\tunip$ is the desired
topological $F$-Jordan decomposition.
\end{proof}

Now we show that we can drop the tameness hypotheses of
Proposition \ref{prop:x-depth} and Lemma
\ref{lem:top-Jordan-to-Jordan}.

\begin{cor}
\label{cor:abs-F-ss-tame}
If $\gamma \in G$ is absolutely $F$-semisimple, then it is
$F$-tame.
\end{cor}

\begin{proof}
By Remark \ref{rem:top-F-ss-unip-bdd}, $\gamma$ is bounded.
By Proposition \ref{prop:tJd=cpct}, there is a topological
$F$-Jordan decomposition
$\gamma = \gamma\tsemi\gamma\tunip$
with $\gamma\tsemi$ $F$-tame.
By Proposition \ref{prop:unique-top-F-Jordan},
$\gamma = \gamma\tsemi$.
\end{proof}

The following rather technical result,
which is now an immediate consequence of
Lemma \ref{lem:top-Jordan-to-Jordan},
Lemma \ref{lem:fixed-pts},
and Proposition \ref{prop:tJd=cpct},
is really the heart of
the paper.  It should be viewed as a quite precise
existence result about topological $F$-Jordan
decompositions.

\begin{thm}
\label{thm:top-F-Jordan-projects}
If $x \in \BB(\bG, F)$ and $\gamma \in \stab_G(x)$,
then there is a topological $F$-Jordan decomposition
$\gamma = \gamma\tsemi\gamma\tunip$ such that $\gamma\tsemi$
and $\gamma\tunip$ project to the semisimple and unipotent
parts, respectively, of the image of $\gamma$ in $\ms G_x(\ff)$.
\end{thm}

\begin{rk}
If $F$ is an algebraic extension of a locally compact field,
then the proof of Theorem \ref{thm:top-F-Jordan-projects}
can be considerably simplified.
Indeed, in this case $G$ is ind-locally-pro-$p$, by Remark
\ref{rem:ILPF}; so,
by Propositions \ref{prop:unique-top-p-Jordan_closure} and
\ref{prop:intrinsic-tu}, an element $\gamma \in \stab_G(x)$
has a topological $p$-, hence $F$-, Jordan decomposition
$\gamma =\gamma\tsemi\gamma\tunip$
with $\gamma\tsemi, \gamma\tunip \in \stab_G(x)$.
Now Lemma \ref{lem:top-Jordan-to-Jordan} shows that the
images of $\gamma\tsemi$ and $\gamma\tunip$ in
$\ms G_x(\ff)$ are as desired.
\end{rk}

\begin{cor}
\label{cor:pro-unip-rad}
For $x \in \BB(\bG, F)$, any normal subgroup of $G_x$
consisting entirely of topologically $F$-unipotent elements lies
in $G_x^+$\,.
\end{cor}

\begin{proof}
Suppose that $H \subseteq G_x$ is normal and consists
entirely of topologically $F$-unipotent elements.
By Lemma \ref{lem:top-Jordan-to-Jordan},
the image of $H$ in $\ms G_x\conn(\ff)$
consists entirely of unipotent elements.
Denote by \ms H its Zariski closure in $\ms G_x\conn$.
Then $\ms H\conn$ is a connected, normal, unipotent subgroup of
the reductive group $\ms G_x\conn$, hence trivial.
By Lemma 22.1 of
\cite{borel:linear}, \ms H is central in $\ms G_x\conn$,
hence consists entirely of semisimple elements.
Since we have already observed that it consists entirely of
unipotent elements, \ms H is trivial.
\end{proof}

We already have an existence result (Proposition
\ref{prop:tJd=cpct}) for topological
$F$-Jordan decompositions modulo the trivial group.
The next result handles such decompositions modulo any group
\mo N, for some fields $F$.

\begin{pn}
\label{prop:tJdZ=cpct-mod-Z}
Suppose that $F$ is an algebraic extension of a locally
compact field.
Then the following statements about an element
$\gamma \in G$ are equivalent.
\begin{enumerate}
\item\label{prop:tJdZ=cpct-mod-Z_p} $\gamma$ has a
topological $p$-Jordan decomposition modulo $N$.
\item\label{prop:tJdZ=cpct-mod-Z_F} $\gamma$ has a
topological $F$-Jordan decomposition modulo \mo N.
\item\label{prop:tJdZ=cpct-mod-Z_bdd} $\gamma$ is bounded
modulo \mo N.
\end{enumerate}
\end{pn}

\begin{proof}
By Corollary \ref{cor:abs-ss-loc-cpct} and Lemma
\ref{lem:top-Jordan}, a topological $p$-Jordan decomposition
modulo $N$ is a topological $F$-Jordan decomposition
modulo \mo N, and conversely.
Thus the equivalence (\ref{prop:tJdZ=cpct-mod-Z_p})
$\Leftrightarrow$ (\ref{prop:tJdZ=cpct-mod-Z_F}) is clear.

Denote by $g \mapsto \ol g$ the natural map $\bG \to \tbG$.
By Remark \ref{rem:ILPF}, \tG is ind-locally-pro-$p$, so
we have by Proposition \ref{prop:intrinsic-tu} that \ogamma has a
topological $p$-Jordan decomposition if and only if \ogamma is
compact; that is, if and only if $\gamma$ is compact modulo $N$
(equivalently, by Remark \ref{rem:bdd=cpct}, bounded modulo \mo N).
Thus, to prove the equivalence (\ref{prop:tJdZ=cpct-mod-Z_p})
$\Leftrightarrow$ (\ref{prop:tJdZ=cpct-mod-Z_bdd}),
it suffices to prove that $\gamma$ has a
topological $p$-Jordan decomposition modulo $N$ if and only if
\ogamma has a topological $p$-Jordan decomposition.

The `only if' direction is easy.  For the `if' direction,
suppose that \ogamma has a topological $p$-Jordan decomposition
$\ogamma = (\ogamma)\tsemi(\ogamma)\tunip$.
Denote by $H$ the closure of the group generated by
$\gamma$.
By Remark \xref{exp-rem:G-image} of
\cite{adler-spice:good-expansions}, the image in \tG of
$H$ is closed.  By Proposition \ref{prop:unique-top-p-Jordan_closure},
$(\ogamma)\tsemi$ belongs to this image.  Let
$\gamma\tsemi$ be a preimage of $(\ogamma)\tsemi$ in $H$.
Then $\gamma\tunip := \gamma\tsemi\inv\gamma$ is a preimage
of $(\ogamma)\tunip$, and clearly $\gamma\tsemi$ and
$\gamma\tunip$ commute.
Thus, $\gamma = \gamma\tsemi\gamma\tunip$ is a topological
$p$-Jordan decomposition modulo $N$.
\end{proof}

We close by showing that the ``common torus'' condition of
Definition \ref{defn:top-F-Jordan} can be omitted.

\begin{pn}
\label{prop:top-F-Jordan-defn}
Suppose that $\gamma \in G$
and
$(\gamma\tsemi, \gamma\tunip)$ is a pair of commuting elements
of $G$ such that
\begin{itemize}
\item $\gamma = \gamma\tsemi\gamma\tunip$,
\item $\gamma\tsemi$ is absolutely $F$-semisimple modulo \mo N,
and
\item $\gamma\tunip$ is topologically $F$-unipotent modulo
\mo N.
\end{itemize}
Then $(\gamma\tsemi, \gamma\tunip)$ is a topological
$F$-Jordan decomposition modulo \mo N.
\end{pn}

\begin{proof}
It remains only to show that the images of $\gamma\semi$ and
$\gamma\tsemi$ in $\tbG(\ol F)$ belong to a common
\ol F-torus there.
We will show the equivalent statement that the images of
$\gamma\tsemi$ and $(\gamma\tunip)\semi$ belong to a common
\ol F-torus.
Clearly, it suffices to assume that \mo N is the trivial
subgroup, so we do so.
By Remark \ref{rem:top-F-ss-unip-ascent}, the hypotheses
remain valid if we replace $F$ by a finite extension,
so we may, and hence do, make such
replacements as necessary.

Upon replacing $F$ by a finite extension, we may, and hence
do, assume that $\gamma\semi$ (hence $(\gamma\tunip)\semi$)
lies in $G$
and $\gamma\tsemi$ is $F$-split.
By Lemma \ref{lem:top-F-unip-Jordan},
$(\gamma\tunip)\semi$ is topologically $F$-unipotent.
Thus, upon replacing $F$ by a finite
extension, we may, and hence do, assume that
$(\gamma\tunip)\semi \in G_{0+}$\,.
By Remark \ref{rem:top-F-ss-unip-bdd} and
Lemma \ref{lem:bounded-and-top-F-unip}, there is an
element $x \in \BB(\bG, F)$ such that
$\gamma\tsemi\dota x = x$
and $\gamma\tunip \in G_x^+$\,.
By Proposition \xref{exp-prop:compatibly-filtered-tame-rank}
of \cite{adler-spice:good-expansions},
$\bH := C_\bG(\gamma\tsemi)$ is a compatibly filtered
$F$-subgroup of \bG, in the sense of Definition
\xref{exp-defn:compatibly-filtered} of \emph{loc.\ cit.}
Thus, $\BB(\bH, F)$ may be regarded (non-canonically) as a
subset of $\BB(\bG, F)$ in such a way that
$G_z^+ \cap H = H_z^+$ for $z \in \BB(\bH, F)$.
By Proposition \ref{prop:x-depth}, we have that
$x \in \BB(\bH, F)$,
so
$(\gamma\tunip)\semi \in G_x^+ \cap H = H_x^+
	\subseteq H\conn$\,.
In particular, $(\gamma\tunip)\semi$ belongs to some maximal
$F$-torus \bT in \bH.
Since $\gamma\tsemi$ is central in \bH, it also belongs to
\bT.
\end{proof}

\subsection{Topological unipotence and tameness}

We have already seen that an absolutely $F$-semisimple
element is $F$-tame (see Corollary \ref{cor:abs-F-ss-tame}).
Of course, a topologically $F$-unipotent element need not be
$F$-tame (or even semisimple).
In the next result, we see that, for $F$-tame topologically
$F$-unipotent elements, it is not necessary to introduce the
finite extension $E/F$ of Definition
\ref{defn:top-F-ss-unip}.

\begin{prop}
\label{prop:tame-tunip}
The topologically $F$-unipotent part of a bounded and $F$-tame
element belongs to $G_{0+}$\,.
\end{prop}

\begin{proof}
Let $\gamma$ be a bounded and $F$-tame element.
By Proposition \ref{prop:tJd=cpct}, it has a topological $F$-Jordan
decomposition
$\gamma = \gamma\tsemi\gamma\tunip$
with $\gamma\tsemi$ $F$-tame.
By Lemma \xref{exp-lem:domain-field-ascent} of
\cite{adler-spice:good-expansions}, we may, and hence do,
replace $F$ by a finite tame extension
so that \bG is $F$-quasisplit, and
$\gamma$ and $\gamma\tsemi$ belongs to $F$-split tori.

If $p = 0$, then let $E/F$ be a finite
extension such that $\gamma\tunip \in \bG(E)_{0+}$\,.
By another application of Lemma
\xref{exp-lem:domain-field-ascent} of \emph{loc.\ cit.},
$\gamma\tunip \in G_{0+}$\,.

If $p > 0$, then let
\begin{itemize}
\item \bS be a maximal $F$-split torus containing $\gamma$,
\item \bT the maximal torus containing \bS,
and
\item $E/F$ the splitting field of \bT.
\end{itemize}
By Lemma \ref{lem:top-F-Jordan-central},
$\gamma\tsemi \in Z(C_G(\gamma)) \subseteq T$.
Thus, \bS commutes with $\gamma\tsemi$,
hence is a maximal $F$-split torus in
$C_\bG(\gamma\tsemi)\conn$.
By Lemma \xref{exp-lem:split-in-center} of
\cite{adler-spice:good-expansions}, we have $\gamma\tsemi \in S$.
Therefore $\gamma\tunip \in S$ also.
Since $\gamma\tunip$ is bounded (by Remark
\ref{rem:top-F-ss-unip-bdd}) and
$S\subb = S_0$ (by Remark \ref{rem:torus-filt}),
we have $\gamma\tunip \in S_0$\,.
By Lemma \ref{lem:top-unip}, $\gamma\tunip$ is topologically
$p$-unipotent, hence has $p$-power order modulo
$S_{0+}$\,.
On the other hand,
$\ms S(\ff) = S_0/S_{0+}$ is the group
of \ff-rational points
of an \ff-split torus, hence contains no non-trivial
elements of $p$-power order.
That is,
$\gamma\tunip \in S_{0+}$\,.
By Lemma 2.4 of \cite{rapoport:T1-is-T0}, we have
$S_0 \subseteq T_0$\,.
By Remark \ref{rem:torus-filt}, we have
$\bS(E)_{0+} \subseteq \bT(E)_{0+}$\,.
By Lemmata \xref{exp-lem:torus-field-descent}
and \xref{exp-lem:levi-descent} of
\cite{adler-spice:good-expansions},
$$
S_{0+}
= S_0 \cap \bS(E)_{0+}
\subseteq T_0 \cap \bT(E)_{0+}
= T_{0+}
\subseteq G_{0+}\,.
$$
In particular, $\gamma\tunip \in G_{0+}$\,, as desired.
\end{proof}

\subsection{Lifting}

In this subsection, put $\tbG = \bG/Z(\bG)\conn$.
(This is consistent with the notation in the earlier part of
the paper, as long as we take $\mo N = Z(\bG)\conn$.)
Denote by $g \mapsto \ol g$ the natural map $\bG \to \tbG$.

We show that elements of \tG which are absolutely
$F$-semisimple or topologically $F$-unipotent can be lifted,
upon passing to suitable finite extensions $E/F$, to elements
of $\bG(E)$ which are absolutely $E$-semisimple or
topologically $E$-unipotent, respectively.

\begin{pn}
\label{prop:ts-lift}
If $\gamma \in G$ is absolutely $F$-semisimple modulo
$Z(\bG)\conn$,
then there is a finite separable extension $E/F$
such that $\gamma Z(\bG)\conn(E)$
contains an absolutely $E$-semisimple element.
\end{pn}

\begin{proof}
Since \ogamma is semisimple,
the unipotent part of $\gamma$ lies in $Z(\bG)\conn(\ol F)$,
hence is trivial.  That is, $\gamma$ is semisimple.

Let
\begin{itemize}
\item \bT be an $F$-torus such that $\gamma \in T$,
\item $E/F$ the splitting field of \bT,
and
\item $e_\gamma$ the homomorphism
$\bX^*(\bT) \to \ord(E\cross)$ sending
$\chi \in \bX^*(\bT)$ to $\ord(\chi(\gamma))$.
\end{itemize}
By Remark \ref{rem:top-F-ss-unip-bdd},
$\gamma$ is bounded modulo $Z(\bG)\conn$, so
$e_\gamma$ is trivial on
$\wtilde\bY^* := \bX^*(\bT/(Z(\bG)\conn \cap \bT))$.
Since
$\ord(E\cross)$ is torsion-free and
$\bY^* := \bX^*(\bT/(Z(\bG) \cap \bT)\conn)$ has finite index
in $\wtilde\bY^*$,
also $e_\gamma$ is trivial on $\bY^*$,
hence induces a homomorphism
$\lambda$
from $\bX^*(\bT)/\bY^* \cong \bX^*((Z(\bG) \cap \bT)\conn)$
to $\ord(E\cross)$.
By choosing a uniformiser for $E$, hence an isomorphism
$\ord(E\cross) \cong \Z$, we may, and hence do, regard
$\lambda$ as an element of
$\bX_*((Z(\bG) \cap \bT)\conn)$.
Denote by $z \in Z(\bG)\conn$ the value of $\lambda$ at the chosen
uniformiser, so that
$\ord(\chi(z)) = e_\gamma(\chi) = \ord(\chi(\gamma))$ for all
$\chi \in \bX^*(\bT)$.
Then $\delta := \gamma z\inv \in \bT(E)$ is bounded.

By Proposition \ref{prop:tJd=cpct},
there exists a topological $E$-Jordan decomposition
$\delta = \delta\tsemi\delta\tunip$.
Notice that
$\ol{\delta\tsemi}$ and $\ol{\delta\tunip}$
belong to a common \ol F-torus
(namely, the image in \tbG of any \ol F-torus in
\bG containing both $\delta\tsemi$ and $\delta\tunip$).
Clearly, $\ol{\delta\tsemi}$ is absolutely
$E$-semisimple.
As in the proof of Proposition
\ref{prop:unique-top-F-Jordan}, the character values of
$\delta\tunip$ lie in $K\cross_{0+}$ for some
finite extension $K/E$.
The character values of $\ol{\delta\tunip}$, being a subset of those
of $\delta\tunip$, thus also belong to $K\cross_{0+}$\,.
By replacing $K$ by a further finite (separable) extension
if necessary, we may, and hence do, assume
that $\bT/Z(\bG)\conn$ is $K$-split,
so that Remark \ref{rem:torus-filt} gives
$\ol{\delta\tunip} \in (\bT/Z(\bG)\conn)(K)_{0+}$\,.
By Lemma \xref{exp-lem:levi-descent}
of \cite{adler-spice:good-expansions}, we have that
$\ol{\delta\tunip} \in \tbG(K)_{0+}$\,,
so $\ol{\delta\tunip}$ is topologically $E$-unipotent.
That is,
$\ogamma = \ol\delta =
\ol{\delta\tsemi}\dotm\ol{\delta\tunip}$
is a topological $E$-Jordan decomposition of \ogamma.
By Proposition \ref{prop:unique-top-F-Jordan},
$\ol{\delta\tunip} = \ol 1$, so
$\delta\tunip \in Z(\bG)\conn(E)$.
Thus $\gamma Z(\bG)\conn(E)$ contains an absolutely
$E$-semisimple element, namely
$\delta\tsemi = \gamma z\inv\delta\tunip\inv$, as desired.
\end{proof}

\begin{rk}
\label{rem:tame-ts-lift}
The field $E/F$ occurring in Proposition \ref{prop:ts-lift}
may be taken to be the splitting field for any $F$-torus
containing $\gamma$.  In particular, if $\gamma$ is
$F$-tame, then $E/F$ may be chosen to be tame.
(Notice that Corollary \ref{cor:abs-F-ss-tame} only guarantees
that \ogamma, not $\gamma$ itself, is $F$-tame.)
We do not know an equally satisfactory answer to when the
field extension $E/F$ in the next proposition may be taken
to be tame.
\end{rk}

\begin{pn}
\label{prop:tu-lift}
If $\gamma \in G$ is topologically $F$-unipotent modulo
$Z(\bG)\conn$,
then there is a finite extension $E/F$
such that $\gamma Z(\bG)\conn(E)$
contains a topologically $E$-unipotent element.
\end{pn}

\begin{proof}
By Remark \ref{rem:top-F-ss-unip-ascent}, we may, and hence
do, replace $F$ by a finite extension
so that \bG is $F$-split and $\ogamma \in \tG_{0+}$\,;
say $x \in \BB(\bG, F)$ is such that
$\gamma \in \tG_\ox^+$\,,
and \bT is an $F$-split maximal torus in \bG whose apartment
contains $x$.

It suffices to show that the
image of $G_x^+$ under the natural map
$G \to \tG$ includes $\tG_\ox^+$\,.
By Remark \xref{exp-rem:actually-product} of
\cite{adler-spice:good-expansions},
since the affine root subgroups of \tbG are naturally isomorphic to
those of \bG, it suffices to show that
the image of $T_{0+}$ includes
$\wtilde T_{0+}$ (where $\wtilde\bT := \bT/Z(\bG)\conn$).
The following square commutes:
$$
\begin{CD}
T_{0+} @>>> \Hom_\Z(\bX^*(\bT), F\cross_{0+}) \\
@VVV @VVV \\
\wtilde T_{0+} @>>> \Hom_\Z(\bX^*(\wtilde\bT), F\cross_{0+})
\end{CD}
$$
(where the vertical maps are the obvious ones, the
top horizontal map takes $t \in T_{0+}$ to the
``evaluation at $t$'' homomorphism, and the
bottom horizontal map is the analogous map for $\wtilde\bT$).
By Remark \ref{rem:torus-filt}, the
top and bottom horizontal arrows are isomorphisms.
The cokernel of the right-hand vertical map is
$\Ext^1_\Z(\bX^*(Z(\bG)\conn), F\cross_{0+})$, which is
trivial since $\bX^*(Z(\bG)\conn)$ is a free \Z-module.
Thus the left-hand map is surjective, as desired.
\end{proof}

\end{document}